\documentclass[letterpaper, 10 pt, conference]{ieeeconf}%
\usepackage{graphics}
\usepackage{epsfig}
\usepackage{cite}
\usepackage{breakurl}
\usepackage{amsmath}
\usepackage{times}

\usepackage{amsthm}
\usepackage{bm}
\usepackage{amssymb}
\usepackage{amsfonts}
\usepackage{bm}
\usepackage{amssymb}
\usepackage{amsfonts}
\usepackage{graphicx}
\usepackage{makecell,multirow,diagbox}
\usepackage{footnote}
\usepackage[normalem]{ulem}
\usepackage{color}
\usepackage{dsfont}
\usepackage[hyphens]{url}
\usepackage[hidelinks]{hyperref}
\usepackage[ruled,linesnumbered]{algorithm2e}
\usepackage[dvipsnames]{xcolor}
\let\labelindent\relax
\usepackage{enumitem}%
\usepackage{booktabs, tabularx}%
\providecommand{\U}[1]{\protect\rule{.1in}{.1in}}
\providecommand{\U}[1]{\protect\rule{.1in}{.1in}}
\IEEEoverridecommandlockouts
\newtheorem{customthm}{Theorem}
\newtheorem{customlem}{Lemma}
\newtheorem{customass}{Assumption}

\newtheorem{customdef}{Definition}
\newtheorem{customrem}{Remark}
\allowdisplaybreaks[0]
\hypersetup{breaklinks=true}
\begin{document}

\title{{\LARGE \textbf{Decentralized Optimal Control of Connected Automated Vehicles
at Signal-Free Intersections Including Comfort-Constrained Turns and Safety
Guarantees}}}
\author{Yue Zhang, and Christos G. Cassandras \thanks{Supported in part by NSF under
grants ECCS-1509084 and CNS-1645681, by AFOSR under grant FA9550-15-1-0471, by
ARPA-E's NEXTCAR program under grant DE-AR0000796 and by the MathWorks.}
\thanks{Y. Zhang and C.G. Cassandras are with the Division of Systems
Engineering and Center for Information and Systems Engineering, Boston
University, Boston, MA 02215 USA (e-mail: joycez@bu.edu; cgc@bu.edu).} }
\maketitle

\begin{abstract}
We extend earlier work for optimally controlling Connected Automated Vehicles
(CAVs) crossing a signal-free intersection by including all possible turns
taken so as to optimize a passenger comfort metric along with energy and
travel time minimization. We show that it is possible to achieve this goal in
a decentralized manner with each CAV solving an optimal control problem, and
derive explicit solutions that guarantee collision avoidance and safe distance
constraints within a control zone. We investigate the associated tradeoffs
between minimizing energy and vehicle travel time, as well as the passenger
comfort metric and include extensive simulations to illustrate this framework.

\end{abstract}



\section{Introduction}

\label{intro}

Intersections in transportation systems are one of the major traffic control
challenges as they account for a large part of accidents and of overall road
congestion. As documented in \cite{scorecard2015texas}, in 2014, traffic
congestion caused vehicles in urban areas to spend 6.9 billion extra hours on
the road at a cost of an extra 3.1 billion gallons of fuel, resulting in a
total cost estimated at \$160 billion. From a control and optimization
standpoint, the goal is to develop efficient traffic management methods so as
to reduce congestion and increase safety with minimal impact on the existing
infrastructure. This is typically accomplished through tighter spacing of
vehicles \cite{Malikopoulos2013}, \cite{Margiotta2011} which can alleviate
congestion, reduce energy use and emissions, and improve safety under proper
control. Forming \textquotedblleft platoons\textquotedblright\ of vehicles is
a popular system-level approach that gained momentum in the 1990s
\cite{Shladover1991,Rajamani2000}. More recently, a study in
\cite{tachet2016revisiting} indicated that transitioning from intersections
with traffic lights to autonomous ones has the potential of doubling capacity
and reducing delays.

To date, traffic light control is the prevailing method for coordinating
conflicting traffic flows and ensure road safety in urban areas. Recent
technological developments include designing adaptive traffic light control
systems that can dynamically adjust the signal timing to various context,
e.g., \cite{fleck2016adaptive} and references therein. However, in addition to the
obvious infrastructure cost of traffic lights, the efficiency and safety
offered by such signaling methods is limited, thus motivating research efforts
to explore new approaches capable of enabling smoother traffic flow while
ensuring safety.

Connected Automated Vehicles (CAVs), also referred to as autonomous or
self-driving vehicles, provide the most intriguing opportunity for better
traffic conditions in a transportation network and for improving traffic flow.
CAVs can assist drivers in making better operating decisions or they can do so
in a fully automated way so as to improve safety and reduce pollution, energy
consumption, and travel delays. One of the very early efforts in this
direction was proposed in \cite{Athans1969} and \cite{Levine1966} where the
merging problem was formulated as a linear optimal regulator to control a
single string of vehicles. The key features of an automated intelligent
vehicle-highway system (IVHS) and a related control system architecture are
discussed in \cite{varaiya1993smart}. More recently, several research efforts
have been reported in the literature for CAV coordination at intersections.
Dresner and Stone \cite{Dresner2004} proposed a reservation-based scheme for
centralized automated vehicle intersection management while no turns are
allowed. Since then, numerous research efforts have explored safe and
efficient control strategies, e.g., \cite{Dresner2008, DeLaFortelle2010,
Huang2012}. Several efforts have focused on minimizing vehicle travel time
under collision-avoidance constraints \cite{Li2006,Yan2009,Zhu2015,Zohdy2012,
Lee2012, Miculescu2014}. Zohdy et al. \cite{Zohdy2012} presented an approach
based on Cooperative Adaptive Cruise Control (CACC) for minimizing
intersection delay and hence maximizing the throughput. Lee and Park
\cite{Lee2012} considered minimizing the overlap between vehicle positions.
Miculescu and Karaman \cite{Miculescu2014} have studied intersections as
polling systems and determined a sequence of times assigned to vehicles on
each road which provides provable guarantees on safety and expected wait time.
Reducing energy consumption is another desired objective which has been
considered in recent literature \cite{gilbert1976vehicle, hooker1988optimal,
hellstrom2010design, li2012minimum}. Hellstrom \cite{hellstrom2010design}
proposed an energy-optimal control algorithm for heavy diesel trucks by
utilizing road topography information. Based on Vehicle-to-Vehicle (V2V)
communication, a minimum energy control strategy is investigated in
car-following scenarios in \cite{li2012minimum}. A detailed discussion of
recent advances in this area can be found in \cite{rios2017survey}.

Our earlier work \cite{ZhangMalikopoulosCassandras2016} and
\cite{malikopoulos2018decentralized} established a decentralized optimal control framework
for coordinating on line a continuous flow of CAVs crossing an urban
intersection without using explicit traffic signaling, assuming no left and
right turns are allowed. For each CAV, an energy minimization optimal control
problem is formulated where the time to cross the intersection is first
determined through a throughput maximization problem. We solved this problem
for each CAV entering a specified \textit{Control Zone} (CZ) so that the
acceleration/deceleration of the CAV is controlled until it reaches a
\textit{Merging Zone} (MZ) where the potential of lateral collisions exists.
To ensure safety, we required CAVs to have a constant speed through the MZ
while also maintaining a safe distance between them to avoid rear-end
collisions. The presence of hard safety constraints makes it challenging to
ensure the existence of a feasible solution to each such decentralized
problem. Therefore, we also established in
\cite{ZhangMalikopoulosCassandras2016} the conditions under which such
feasible solutions exist and showed that they can be enforced through an
appropriately designed Feasibility Enforcement Zone (FEZ) that precedes the CZ
in \cite{Zhang2016}.

The inclusion of left and right turns in this framework creates significant
complications. Aside from the added complexity of ensuring collision
avoidance, we can no longer expect each CAV to maintain a constant speed in
the MZ. Instead, we must allow a CAV to vary its speed depending on the turn.
In addition, passenger \emph{comfort} in taking such turns becomes an
essential component of an automated trajectory design. The problem of
coordinating CAVs at intersections including left and right turns has been
addressed by Kim and Kumar \cite{kim2014mpc} based on an approach using Model
Predictive Control (MPC) to achieve system-wide safety and liveness of
intersection-crossing traffic. Dresner and Stone in
\cite{dresner2005multiagent} considered scenarios that allowed left and right
turns using a reservation-based scheme together with a communication protocol
which may involve a stop sign and traffic lights.

The contribution of this paper consists of extending the optimal control
framework in \cite{ZhangMalikopoulosCassandras2016}. First, unlike the
approach in \cite{ZhangMalikopoulosCassandras2016} where we first solved a
throughput maximization problem to determine a CAV's travel time and then an
energy minimization problem for each CAV, here we formulate a problem in which
a CAV seeks to \emph{jointly minimize} both its travel time through the CZ and
its energy consumption. This allows us to readily quantify the tradeoff
between these two criteria. Second, we formulate and solve a subsequent
optimal control problem for each CAV crossing the MZ in which we allow it to
vary its speed inside the MZ and include a passenger comfort metric when a
turn is taken. The ability to decentralize the optimal control in
\cite{ZhangMalikopoulosCassandras2016} rests on showing that each CAV needs
only information from the CAV which is physically ahead of it and the one
immediately preceding it in entering the CZ. The presence of turns complicates
this simple coordination structure; however, we are still able to show that
each CAV needs only information from a small set of well-defined CAVs among
those preceding it (but not necessarily immediately preceding it) in entering
the CZ. We also derive explicit solutions for the two optimal control
problems, through the CZ and then the MZ, including the possibility of safety,
state and control constraints becoming active. Our analysis includes the
derivation of properties characterizing an optimal control solution, such as
the continuity of the optimal control when an unconstrained optimal trajectory
arc is followed by one with an active constraint, and it allows us to
determine whether an optimal control solution for each CAV is feasible at the
time it enters the CZ.

The paper is structured as follows. In Section II, we review the model in
\cite{ZhangMalikopoulosCassandras2016} and extend it to include left and right
turns. In Section III, we derive the conditions that guarantee safety for each
CAV in terms of its time to reach the MZ constrained by those of CAVs
preceding it in the CZ. In Section IV, we formulate a decentralized optimal
control problem for each CAV that jointly minimizes its travel time and energy
consumption, prove structural properties of optimal trajectories, and derive
an explicit solution for it. In Section V, we formulate and solve another optimization
problem with the objective of jointly minimizing a passenger comfort metric
inside the MZ and its energy consumption. Simulation results are given in
Section VI where the time-energy tradeoff is illustrated, as well as the
comfort-energy tradeoff.


\section{The Intersection Model}

\label{model}

\begin{figure}[ptb]
\centering
\includegraphics[width = 0.7 \columnwidth]{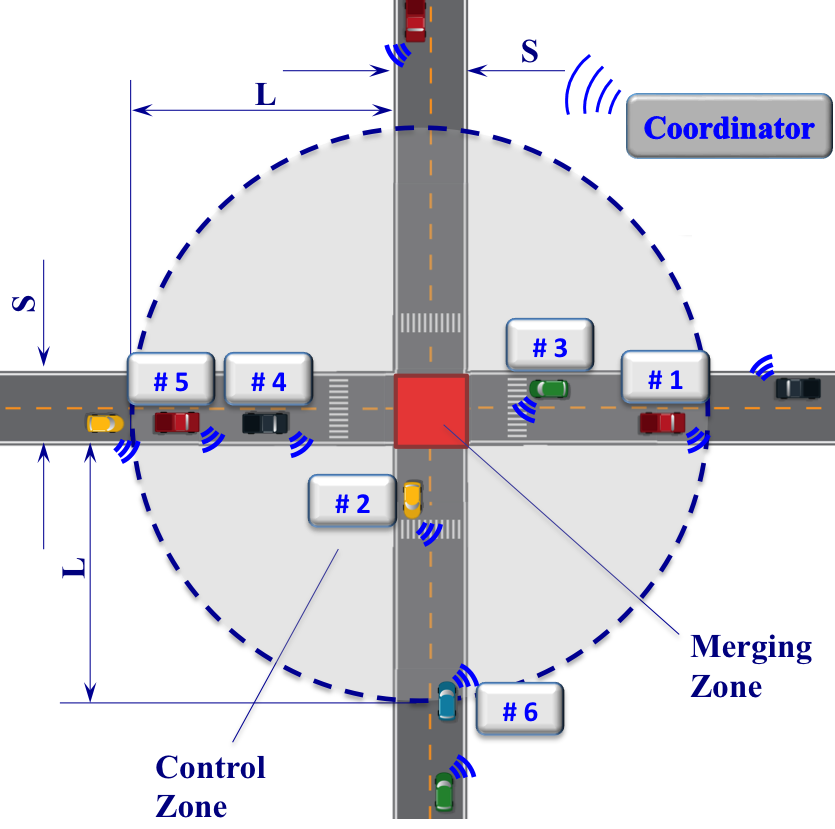}
\caption{Connected Automated Vehicles crossing an urban intersection.}
\label{fig:intersection}
\end{figure}

We begin with a brief review of the model introduced in \cite{ZhangMalikopoulosCassandras2016} and fully developed in
\cite{malikopoulos2018decentralized}. We consider an intersection (Fig.
\ref{fig:intersection}) where the region at its center, assumed to be a square of side $S$, is called \emph{Merging Zone} (MZ) and defines the area of potential lateral CAV collisions. The intersection has a \emph{Control Zone} (CZ) and the road segment from the CZ entry to the CZ exit (i.e., the MZ entry) is referred to as a CZ segment. The length of each CZ segment is $L>S$ and it is assumed to be the same for all entry points to a given CZ. Extensions to asymmetric CZ segments are possible and considered in \cite{Zhang2018sequence}.

We assume the existence of a ``coordinator"  whose task is to handle the information exchanges between CAVs, while each
CAV maintains its own control autonomy. Let $N(t)\in\mathbb{N}$ be the
cumulative number of CAVs which have entered the CZ by time $t$ and formed a
queue that designates the crossing sequence in which these CAVs will enter the
MZ. There is a number of ways to manage such a queue. In
\cite{ZhangMalikopoulosCassandras2016} a strict First-In-First-Out (FIFO)
crossing sequence is assumed, that is, when a CAV reaches the CZ, the
coordinator assigns it an integer value $i=N(t)+1$. This can be relaxed as in
\cite{Zhang2018sequence} to allow for dynamically resequencing CAVs as each
new one arrives, hence maximizing throughput. If two or more CAVs enter a CZ
at the same time, then the corresponding coordinator selects randomly the
first one to be assigned the value $N(t)+1$.

For simplicity, we assume that each CAV is governed by second order dynamics:%
\begin{equation}
\dot{p}_{i}=v_{i}(t)\text{, }~p_{i}(t_{i}^{0})=0\text{; }~\dot{v}_{i}%
=u_{i}(t)\text{, }v_{i}(t_{i}^{0})\text{ given}\label{eq:model2}%
\end{equation}
where $p_{i}(t)\in\mathcal{P}_{i}$, $v_{i}(t)\in\mathcal{V}_{i}$, and
$u_{i}(t)\in\mathcal{U}_{i}$ denote the position, i.e., travel distance since
the entry of the CZ, speed and acceleration/deceleration (control input) of
each CAV $i$. The sets $\mathcal{P}_{i}$, $\mathcal{V}_{i}$ and $\mathcal{U}%
_{i}$ are complete and totally bounded subsets of $\mathbb{R}$. These dynamics
are in force over an interval $[t_{i}^{0},t_{i}^{f}]$, where $t_{i}^{0}$ and
$t_{i}^{f}$ are the times that the vehicle $i$ enters the CZ and exits the MZ
respectively. To ensure that the control input and vehicle speed are within a
given admissible range, the following constraints are imposed:
\begin{equation}%
\begin{split}
u_{i,min} &  \leq u_{i}(t)\leq u_{i,max},\quad\text{and}\\
0 &  \leq v_{min}\leq v_{i}(t)\leq v_{max},\quad\forall t\in\lbrack t_{i}%
^{0},t_{i}^{f}].
\end{split}
\label{speed_accel constraints}%
\end{equation}
To ensure the absence of any rear-end collision throughout the CZ, we impose
the \emph{rear-end safety} constraint
\begin{equation}
s_{i}(t)=p_{k}(t)-p_{i}(t)\geq\delta,\quad\forall t\in\lbrack t_{i}^{0}%
,t_{i}^{f}]\label{rearend}%
\end{equation}
where $k$ is the CAV physically ahead of $i$ and $\delta$ is the \emph{minimal
safe following distance} allowable. An alternative safety constraint often
used is to maintain a safe headway \cite{Rajamani2000}, i.e., a time gap that
is a function of speed. However, since we consider urban intersections, the
average speed does not exhibit significant variations and we can translate the
allowable headway to a safe inter-vehicle distance (replacing (\ref{rearend})
by a headway constraint is still possible at the expense of complicating the
explicit solutions to the optimal control problems we will consider in the sequel).


In this modeling framework, we assume that $(i)$ CAVs follow the crossing
sequence established by the coordinator and that no overtaking, reversing
directions, or lane-changing are allowed, $(ii)$ each vehicle has proximity
sensors and can observe and/or estimate local information that can be shared
with other vehicles, $(iii)$ the decision of each CAV on whether a turn needs
to be made at the MZ is known upon its entry in the CZ and remains unchanged,
and $(iv)$ for each CAV, the speed constraints in
\eqref{speed_accel constraints} and the rear-end safety constraint in
\eqref{rearend} are not active at $t_{i}^{0}$. If this last assumption is
violated, any optimal control solution is obviously infeasible and we must
resort to control actions that simply attempt to satisfy these constraints as
promptly as possible; alternatively, we may impose a Feasibility Enforcement
Zone (FEZ) that precedes the CZ as described in \cite{Zhang2016}.

Finally, we make the following simplifying assumption:

\begin{customass}
For each CAV $i$ exiting the MZ, its speed remains constant for at least a
distance of length $\delta$. \label{ass:cruise}
\end{customass}

This assumption is reasonable since there is no compelling reason for vehicles
to accelerate or decelerate right after they exit the MZ unless an imminent
safety issue is involved.

\subsection{Modeling Turns \label{section:modelturns}}

\label{modeling_turns} 
\begin{figure}[ptb]
\centering
\includegraphics[width= 0.6\columnwidth]{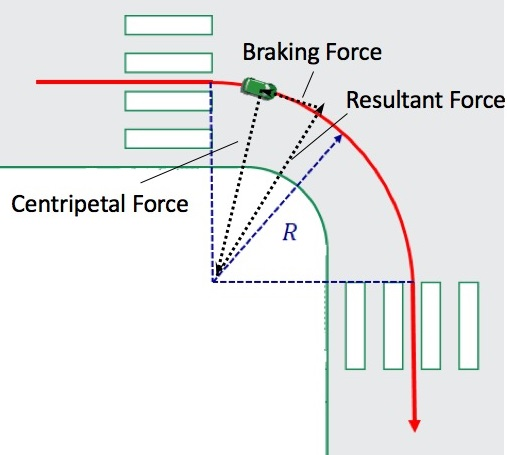} \caption{Vehicle
making a turn.}%
\label{radius}%
\end{figure}

The inclusion of left and right turns needs special
attention in the context of safety as well in ensuring passenger comfort. As a
vehicle turns, it is affected simultaneously by two forces: the braking force
and the centripetal force (see Fig.\ref{radius}). Note that the centripetal
force is provided by the friction force, which depends on the speed of the
vehicle. Hence, to ensure safety, vehicles need to slow down while making
turns. The smaller the turning radius is, the lower the speed needs to be. On
the other hand, to minimize passenger discomfort, vehicles need to maintain a
steady movement, i.e., avoid any abrupt change in the braking force. In other
words, we need to minimize the change in acceleration (i.e., the \emph{jerk}).


Let $d_{i}$ denote the decision of vehicle $i$ on whether a turn is to be made
at the MZ, where $d_{i}=0$ indicates a left turn, $d_{i}=1$ indicates going
straight and $d_{i}=2$ indicates a right turn. Letting $S_{L}$ denote the
length of the left turn trajectory and $S_{R}$ the length of the right turn
trajectory, we assume that $S_{R}<S<S_{L}$ holds. The speed for which an
intersection curve is designed depends on an assigned speed limit, the type of
intersection, and the traffic volume \cite{aashto2001policy}. Generally, the
\textquotedblleft desirable time\textquotedblright\ $\Delta_{i}$ that a
vehicle needs to make a turn at an intersection \cite{aashto2001policy} is
given by%

\begin{equation}
\Delta_{i}=\left\{
\begin{array}
[c]{ll}%
\frac{R_{i}}{\sqrt{15R_{i}(0.01E+F)}}, & \mbox{if $d_i = 0, 2$},\\
\Delta_{i}^{1}, & \mbox{if $d_i = 1$},
\end{array}
\right.  \label{deltat1}%
\end{equation}
where $R_{i}$ is the centerline turning radius (see Fig. \ref{radius}); $E$ is
the super-elevation, which is zero in urban conditions; $F$ is the side
friction factor; and $\Delta_{i}^{1}$ is the time for CAV $i$ going straight.
Therefore, the time $t_{i}^{m}$ when CAV $i$ enters the MZ is directly related
to the time $t_{i}^{f}$ that the vehicle exits the MZ through $\Delta_{i}$:
\begin{equation}
t_{i}^{f}=t_{i}^{m}+\Delta_{i}. \label{deltat2}%
\end{equation}
Note that $\Delta_{i}$ is different for left and right turns since the
associated turning radii $R_{i}$ are different.

\subsection{Relative Location Sets \label{sec:subsets}}

\begin{figure*}[ptb]
\centering
\includegraphics[width= 1.75\columnwidth]{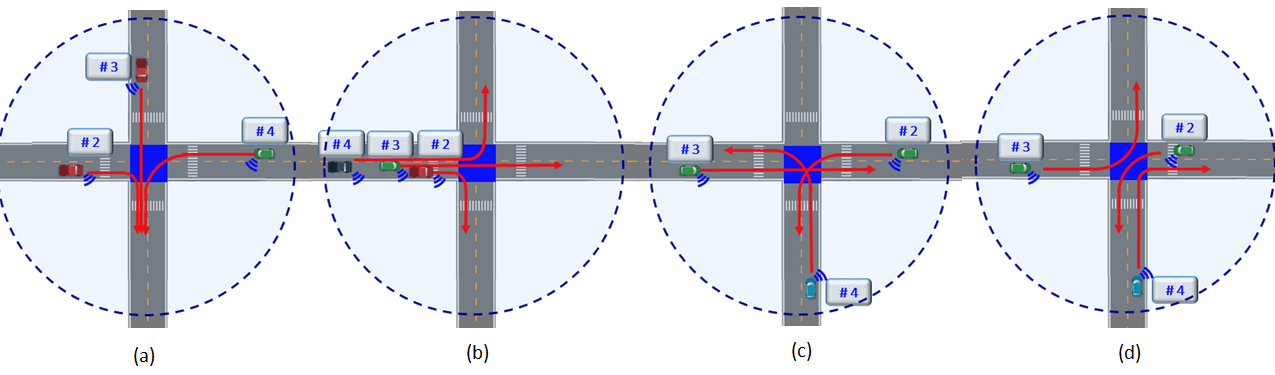}
\caption{Illustration of different subsets of $\mathcal{Q}^{z}(t)$: (a) subset
$\mathcal{E}(t)$; (b) subset $\mathcal{S}(t)$; (c) subset $\mathcal{L}(t)$;
(d) subset $\mathcal{O}(t)$}%
\label{category}%
\end{figure*}

When a CAV enters the CZ and is assigned a unique index $i=N(t)+1$ by the
coordinator, it is also assigned to one and only one of four subsets
$\{\mathcal{E}_{i}(t),$ $\mathcal{S}_{i}(t),$ $\mathcal{L}_{i}(t),$
$\mathcal{O}_{i}(t)\}$ which capture its position relative to CAVs $j<i$.
These subsets are considerably different than those in
\cite{ZhangMalikopoulosCassandras2016} due to the presence of turns and are
defined as follows:

(1) $\mathcal{E}_{i}(t)$ contains all vehicles $j<i$ that can cause a
\emph{rear-end collision} with $i$ at the \textit{end} of the MZ. For example,
in Fig. \ref{category}(a), $\mathcal{E}_{3}(t)$ contains vehicle \#2 as it may
cause a rear-end collision with vehicle \#3 at the end of the MZ, and
$\mathcal{E}_{4}(t)$ contains vehicle \#3 as it may cause a rear-end collision
with vehicle \#4 at the end of the MZ.

(2) $\mathcal{S}_{i}(t)$ contains all vehicles $j<i$ traveling on the same
lane that can cause \emph{rear-end collision} with $i$ at the
\textit{beginning} of the MZ. For example, in Fig. \ref{category}(b),
$\mathcal{S}_{3}(t)$ contains vehicle \#2 as it may cause a rear-end collision
with vehicle \#3 at the beginning of the MZ, and $\mathcal{S}_{4}(t)$ contains
vehicle \#3 as it may cause a rear-end collision with vehicle \#4 at the
beginning of the MZ. Note that if CAVs $j$ and $i$, $j<i$ travel on the same
lane and in the same direction, we have $j\in\mathcal{E}_{i}(t)$; since $j$
belongs to one and only one of these subsets, we cannot have $j\in
\mathcal{S}_{i}(t)$. 

(3) $\mathcal{L}_{i}(t)$ contains all vehicles $j<i$ traveling on different
lanes and towards different lanes that can cause \emph{lateral collision} with
$i$ inside the MZ. For example, in Fig. \ref{category}(c), $\mathcal{L}%
_{3}(t)$ contains vehicle \#2 as it may cause a lateral collision with vehicle
\#3 inside the MZ, and $\mathcal{L}_{4}(t)$ contains vehicle \#3 as it may
cause a lateral collision with vehicle \#4 inside the MZ.

(4) $\mathcal{O}_{i}(t)$ contains all vehicles $j<i$ traveling on different
lanes and towards different directions that cannot cause any lateral collision
with $i$ at the MZ. For example, in Fig. \ref{category}(d), $\mathcal{O}%
_{3}(t)$ contains vehicle \#2 since it cannot cause any collision with vehicle
\#3, and $\mathcal{O}_{4}(t)$ contains vehicle \#3 since it cannot cause any
collision with vehicle \#4.

\subsection{Merging Zone Safety Constraints}

As in \cite{ZhangMalikopoulosCassandras2016}, we consider a First-In-First-Out
(FIFO) queue by imposing the following condition:
\begin{equation}
t_{i}^{f}\geq t_{i-1}^{f},~\text{\ }i>1.\label{eq:fifo}%
\end{equation}
However, as already mentioned, this FIFO structure may be relaxed by
dynamically resequencing CAVs as each new one arrives \cite{Zhang2018sequence}
aiming to maximize throughput and subsequently re-solving each decentralized
CAV problem. As we will see, this is made possible by the relatively modest
computational cost of solving such problems. 

A lateral collision involving CAV $i$ may occur only if some CAV $j\neq i$
belongs to $\mathcal{L}_{i}(t)$. This leads to the following definition:

\begin{customdef}
For each CAV $i\in\mathcal{N}(t)$, the set $\Gamma_{i}$ that includes all time
instants when a lateral collision involving CAV $i$ is possible: $\Gamma
_{i}\triangleq\Big\{t~|~t\in\lbrack t_{i}^{m},t_{i}^{f}]\Big\}$.
\end{customdef}

Consequently, to avoid a lateral collision for any two vehicles $i,j\in
\mathcal{N}(t)$ on different roads, the following constraint should hold:
\begin{equation}
\Gamma_{i}\cap\Gamma_{j}=\varnothing,\text{ \ }\forall t\in\lbrack t_{i}%
^{m},t_{i}^{f}]\text{, \ }j\in\mathcal{L}_{i}(t). \label{eq:lateral}%
\end{equation}

In earlier work \cite{ZhangMalikopoulosCassandras2016}, the speed in the MZ
was considered to be constant, a condition made possible by the fact that, in
the absence of turns, all CAVs have the same trajectory within the MZ. Thus,
$t_{i}^{m}>t_{i-1}^{m}$ could be implicitly ensured by $t_{i}^{f}>t_{i-1}^{f}%
$. Since this is no longer the case, i.e., speeds and position trajectories in
the MZ may be different for CAVs with different driving directions (e.g., CAV
$i$ may turn left while CAV $i-1$ is going straight), the condition $t_{i}%
^{f}\geq t_{i-1}^{f}$ does not imply $t_{i}^{m}\geq t_{i-1}^{m}$. This becomes an issue only when $i$ and $i-1$ are traveling on the same lane.
We shall henceforth reserve the symbol $k$
($k<i$) to denote the index of the CAV physically located ahead of $i$. Then,
for CAVs $i$ and $k$ we impose the following constraint:
\begin{equation}
t_{i}^{m}>t_{k}^{m},~\ i>k\geq1,\label{eq:fifotim}%
\end{equation}
to avoid any rear-end collision at the beginning of the MZ.

\subsection{Terminal Conditions}

\label{terminal_conditions} We now turn our attention to the terminal
conditions that must be enforced on each CAV $i$ at times $t_{i}^{m}$ and
$t_{i}^{f}$ when it enters and exits the MZ respectively. The safety
requirements discussed above imply constraints on these times that depend on
the relative location sets of $i$ as explained next.

(1) Let $e=\max\limits_{j}\{j\in\mathcal{E}_{i}(t)\}$. In this case, CAV $e$
is the vehicle immediately ahead of CAV $i$ in the FIFO queue that may cause a
rear-end collision with $i$ at the end of the MZ. To avoid such a rear-end
collision, $e$ and $i$ should maintain the \textit{minimal safe distance}
$\delta$ defined in (\ref{rearend}) by the time $i$ exits the MZ. Since, by
Assumption \ref{ass:cruise}, each CAV maintains a constant speed after the MZ
exit for at least a distance $\delta$, we set
\begin{equation}
t_{i}^{f}\geq t_{e}^{f}+\frac{\delta}{v_{e}^{f}} \label{Case E}%
\end{equation}
where $t_{e}^{f}$ and $t_{i}^{f}$ are the times that CAVs $e$ and $i$ exits
the MZ respectively, and $v_{e}^{f}$ is the speed of CAV $e$ at the exit of
the MZ. It is easy to see that if the only objective for $i$ is to minimize
its travel time, then it should set $t_{i}^{f}=t_{e}^{f}+\delta/v_{e}^{f}$. In
general, however, as we will see when setting up each CAV's decentralized
optimal control problem in Section IV, (\ref{Case E}) in conjunction with
(\ref{deltat2}) provides only a lower bound constraint on $t_{i}^{m}$ for this
problem:%
\begin{equation}
t_{i}^{m}\geq t_{e}^{f}+\frac{\delta}{v_{e}^{f}}-\Delta_{i} \label{CaseE_tim}%
\end{equation}

Unlike our earlier work \cite{ZhangMalikopoulosCassandras2016} where we
considered a constant speed for all CAVs inside the MZ, the inclusion of left
and right turns forces us to vary the speed $v_{i}(t)$ for $t\in\lbrack
t_{i}^{m},t_{i}^{f}]$. 
The actual values of $v_{i}(t)$, $t\in\lbrack t_{i}%
^{m},t_{i}^{f}]$, are determined by the solution of the optimal control
problem we will define in Section V so as to take a passenger comfort metric
into consideration. For this problem, $t_{i}^{m}$ is the initial time which is
obtained from the decentralized optimal control problem over $[t_{i}^{0}%
,t_{i}^{m}]$ in Section IV.

(2) Let $s=\max\limits_{j}\{j\in\mathcal{S}_{i}(t)\}$. In this case, CAV $s$
is the vehicle immediately ahead of CAV $i$ in the FIFO queue such that it may
cause a rear-end collision at the beginning of the MZ. To guarantee the
rear-end collision constraint does not become active we set%

\begin{equation}
t_{i}^{m}\geq t_{s}^{m}+\Delta_{s}^{\delta}, \label{Case 2 tim}%
\end{equation}
where, in view of (\ref{deltat1}),
\begin{equation}
\Delta_{s}^{\delta}=\left\{
\begin{array}
[c]{ll}%
\frac{\delta}{\sqrt{15R_{i}(0.01E+F)}}, & \mbox{if $d_{s} = 0, 2$},\\
\Delta_{s}^{1\delta}, & \mbox{if $d_{s} = 1$}.
\end{array}
\right.  \label{delta1}%
\end{equation}
is the time vehicle $s$ needs to travel a distance $\delta$ inside the MZ,
where $\Delta_{s}^{1\delta}$ is the time for CAV $s$ traveling a distance
$\delta$ if it goes straight. The time $t_{i}^{f}$ when CAV $i$ will be
exiting the MZ is given by \eqref{deltat2} and, using (\ref{Case 2 tim}), we
get $t_{i}^{f}\geq t_{s}^{m}+\Delta_{s}^{\delta}+\Delta_{i}$. However, if $s$
makes a left turn and $i$ makes a right turn, since $S_{L}>S_{R}$, we may have
$t_{i}^{f}<t_{s}^{f}$, which violates \eqref{eq:fifo}. In that case, we must
set $t_{i}^{f}\geq t_{s}^{f}$. This situation arises whenever the time to
cover the trajectory of $s$ through the MZ is longer that that of $i$.
Combining these two cases, it follows that%

\begin{equation}
t_{i}^{f}\geq\max\{t_{s}^{m}+\Delta_{s}^{\delta}+\Delta_{i},t_{s}^{f}\}
\label{Case S}%
\end{equation}
and, from (\ref{deltat2}), $t_{i}^{m}$ is constrained so that
\[
t_{i}^{m}\geq\max\{t_{s}^{m}+\Delta_{s}^{\delta},t_{s}^{f}-\Delta_{i}\}
\]

(3) Let $l=\max\limits_{j}\{j\in\mathcal{L}_{i}(t)\}$. In this case, CAV $l$
is the vehicle immediately ahead of CAV $i$ in the FIFO queue such that it may
cause a lateral collision inside the MZ. We constrain the MZ to contain only
either CAV $l$ or $i$ so as to avoid such a collision; this is intended to
enhance safety awareness, but it could be modified appropriately, if
necessary. Therefore, CAV $i$ will enter the MZ only after $l$ has exited it,
that is,%

\begin{equation}
t_{i}^{m}\geq t_{l}^{f} \label{lem11}%
\end{equation}
and $t_{i}^{f}$ is obtained through \eqref{deltat2}.

(4) Let $o=\max\limits_{j}\{j\in\mathcal{O}_{i}(t)\}$. In this case, CAV $o$
is the vehicle immediately ahead of CAV $i$ in the FIFO queue such that it
will not generate any collision with $i$ in the MZ, so we only require that
\begin{equation}
t_{i}^{f}\geq t_{o}^{f} \label{lem12}%
\end{equation}
and, from (\ref{deltat2}), $t_{i}^{m}$ is constrained so that
\begin{equation}
t_{i}^{m}\geq t_{o}^{f}-\Delta_{i} \label{lem12_tim}%
\end{equation}

We are now in a position to establish a lower bound on $t_{i}^{m}$, the time
for CAV $i$ to enter the MZ, which will serve as a terminal time constraint in
the decentralized optimal control problem considered in Section IV. In so
doing, we need to also take into account the lower bound $t_{i}^{L}$ imposed
on $t_{i}^{m}$ due to the control and state constraints
\eqref{speed_accel constraints}. To derive this lower bound of $t_{i}^{m}$, we
have to consider two cases, which depend on whether CAV $i$ can reach
$v_{max}$ upon arriving at the MZ: $(i)$ If CAV $i$ enters the CZ at
$t_{i}^{0}$, accelerates with $u_{i,max}$ until it reaches $v_{max}$ and then
cruises at this speed until it leaves the MZ at time $t_{i}^{1}$, it was shown
in \cite{ZhangMalikopoulosCassandras2016} that $t_{i}^{1}=t_{i}^{0}+\frac
{L}{v_{max}}+\frac{(v_{max}-v_{i}^{0})^{2}}{2u_{i,max}v_{max}}$. $(ii)$ If CAV
$i$ accelerates with $u_{i,max}$ but reaches the MZ at $t_{i}^{m}$ with speed
$v_{i}^{m}<v_{max}$, it was shown in \cite{ZhangMalikopoulosCassandras2016}
that $t_{i}^{2}=t_{i}^{0}+\frac{\sqrt{2Lu_{i,max}+(v_{i}^{0})^{2}}-v_{i}^{0}%
}{u_{max}}$. Thus,
\begin{equation}
t_{i}^{L}=t_{i}^{1}\mathds{1}_{v_{i}^{m}=v_{max}}+t_{i}^{2}%
(1-\mathds{1}_{v_{i}^{m}=v_{max}}) \label{t_lw}%
\end{equation}
The following provides a lower bound for $t_{i}^{m}$ ensuring that is feasible.

\begin{customthm}
The lower bound on the time when CAV $i$ can enter the MZ and satisfy all MZ
safety constraints is given by
\begin{equation}
t_{i}^{m}\geq\max\{t_{i}^{L},t_{e}^{f}+\frac{\delta}{v_{e}^{f}}-\Delta
_{i},t_{s}^{m}+\Delta_{s}^{\delta},t_{s}^{f}-\Delta_{i},t_{l}^{f},t_{o}%
^{f}-\Delta_{i}\}, \label{Theorem1}%
\end{equation}
where $e=\max_{j}\{j\in\mathcal{E}_{i}(t)\}$, $s=\max_{j}\{j\in\mathcal{S}%
_{i}(t)\}$, $l=\max_{j}\{j\in\mathcal{L}_{i}(t)\}$, $o=\max_{j}\{j\in
\mathcal{O}_{i}(t)\}$ and $t_{i}^{L}$ is given by (\ref{t_lw}).
\end{customthm}

\emph{Proof. } If $\max\{t_{i}^{L},t_{e}^{f}+\frac{\delta}{v_{e}^{f}}%
-\Delta_{i},t_{s}^{m}+\Delta_{s}^{\delta},t_{s}^{f}-\Delta_{i},t_{l}^{f}%
,t_{o}^{f}-\Delta_{i}\}=t_{i}^{L}$, then $t_{i}^{m}\geq t_{i}^{L}$ ensures
that $t_{i}^{m}$ is feasible since it depends only on the control and state
constraints \eqref{speed_accel constraints}. Let $j<i$ and $j\neq e,s,l,o$.
There are four cases to consider as follows. 

(1) $j\in\mathcal{E}_{i}(t)$. Since $e,j\in\mathcal{E}_{i}(t)$, CAVs $e$ and
$j$ are driving towards the same lane. As $e=\max_{r}\{r\in\mathcal{E}%
_{i}(t)\}$ and $j\neq e$, we have $j<e$ and $j\in\mathcal{E}_{e}(t)$. If
$t_{i}^{m}\geq t_{e}^{f}+\frac{\delta}{v_{e}^{f}}-\Delta_{i}$, then
{\eqref{Case E} holds and we have }$t_{i}^{f}\geq t_{e}^{f}+\frac{\delta
}{v_{e}^{f}}$ and $t_{e}^{f}\geq t_{j}^{f}+\frac{\delta}{v_{j}^{f}}$,
therefore, $t_{i}^{f}>t_{j}^{f}+\frac{\delta}{v_{j}^{f}}$ {which ensures the
absence of a rear-end collision at the end of the MZ and \eqref{Case E}
holds for all }$j\in\mathcal{E}_{i}(t)$. 

(2) $j\in\mathcal{S}_{i}(t)$. Similarly, since we have $s,j\in\mathcal{S}%
_{i}(t)$, CAVs $s$ and $j$ are traveling on the same lane. As $s=\max
_{r}\{r\in\mathcal{S}_{i}(t)\}$ and $j\neq s$, we have $j<s$ and
$j\in\mathcal{S}_{s}(t)$. In this case, a rear-end collision at the beginning
of the MZ is possible. 
If $t_{i}^{m}\geq
\max\{t_{s}^{m}+\Delta_{s}^{\delta},t_{s}^{f}-\Delta_{i}\}$, then
\eqref{Case S} holds and we have $t_{s}^{f}\geq\max\{t_{j}^{m}+\Delta
_{j}^{\delta}+\Delta_{s},t_{j}^{f}\}$ and $t_{i}^{f}\geq\max\{t_{s}^{m}%
+\Delta_{s}^{\delta}+\Delta_{i},t_{s}^{f}\}$ which leads to $t_{i}^{f}\geq
t_{j}^{f}$. Therefore, \eqref{eq:fifo} also holds. 

(3) $j\in\mathcal{L}_{i}(t)$. In this case, a lateral collision is possible
between $i$ and $j$. Given $t_{i}^{m}\geq t_{l}^{f}$, unlike cases (1) and
(2), we cannot determine which subset $j$ belongs to with respect to CAV $l$.
To explore the relationship between CAV $l$ and $j$ explicitly, there are four
subcases to consider as follows. $(i)$ $j\in\mathcal{E}_{l}(t)$. According to
\eqref{Case E}, we have $t_{l}^{f}>t_{j}^{f}$, therefore, $t_{i}^{m}>t_{j}%
^{f}$ and lateral collision avoidance is ensured. $(ii)$ $j\in\mathcal{S}%
_{l}(t)$. According to \eqref{Case S}, we have $t_{l}^{f}\geq\max\{t_{j}%
^{m}+\Delta_{j}^{\delta}+\Delta_{l},t_{j}^{f}\}$ which leads to $t_{l}^{f}\geq
t_{j}^{f}$, hence $t_{i}^{m}\geq t_{j}^{f}$ and lateral collision avoidance is
ensured. $(iii)$ $j\in\mathcal{L}_{l}(t)$. From (\ref{lem11}), $t_{l}^{m}\geq
t_{j}^{f}$ and $t_{i}^{m}\geq t_{l}^{f}$, and it follows that $t_{i}^{m}%
>t_{j}^{f}$ which ensures the absence of a lateral collision. $(iv)$
$j\in\mathcal{O}_{l}(t)$. Based on \eqref{lem12}, we have $t_{l}^{f}\geq
t_{j}^{f}$. Hence, we still have $t_{i}^{m}\geq t_{j}^{f}$ which still ensures
lateral collision avoidance. 

(4) $j\in\mathcal{O}_{i}(t)$. 
According to the definition of $\mathcal{O}_{i}(t)$, $j$ cannot collide with $i$. To summarize, as long
as (\ref{Theorem1}) holds, we can guarantee that all MZ constraints are satisfied.

\hfill$\blacksquare$

Corresponding to the lower bound of terminal time $t_{i}^{L}$, there also
exists the upper bound $t_{i}^{U}$, that is,
\begin{equation}
t_{i}^{U}=t_{i}^{3}\mathds{1}_{v_{i}^{m}=v_{min}}+t_{i}^{4}%
(1-\mathds{1}_{v_{i}^{m}=v_{min}}) \label{t_up}%
\end{equation}
where $v_{i}(t_{i}^{m})=\sqrt{2Lu_{min}+(v_{i}^{0})^{2}},$ and $t_{i}%
^{3}=t_{i}^{0}+\frac{L}{v_{min}}+\frac{(v_{min}-v_{i}^{0})^{2}}{2u_{min}%
v_{min}}$ and $t_{i}^{4}=t_{i}^{0}+\frac{v_{i}(t_{i}^{m})-v_{i}^{0}}{u_{min}}$
are derived in a similar way as $t_{i}^{1}$ and $t_{i}^{2}$ in \eqref{t_lw}
respectively. Based on \eqref{t_up}, the following upper bound constraint
applies:
\begin{equation}
t_{i}^{m}\leq t_{i}^{U} \label{tm_up_constraint}%
\end{equation}

It follows directly from (\ref{Theorem1}) that if each CAV seeks to minimize
its travel time through the intersection it should select a time $t_{i}^{m}$
given by the lower bound on the right-hand-side. The CAV can then formulate a
fixed terminal time optimal control problem with a given objective function
and solve this problem in a decentralized way since the value of $t_{i}^{m}$
depends only on information associated with four CAVs (with indices given by
$e$, $s$, $l$, $o<i)$ whose optimal trajectories have already been determined
prior to the arrival of CAV $i$ at the CZ. This approach was followed in
earlier work \cite{ZhangMalikopoulosCassandras2016} where no turns were
considered. In what follows, we formulate instead a problem that \emph{jointly
aims to minimize the travel time of CAV }$i$\emph{ along with a measure of its
energy consumption}. Therefore, the optimal value of $t_{i}^{m}$ is obtained
as part of this problem's solution and the value of Theorem 1, i.e.,
(\ref{Theorem1}), is to provide a constraint ensuring that this value is feasible.


\section{Optimal Control of CAVs in the CZ \label{CZ}}

The objective of each CAV inside the CZ, i.e., over $[t_{i}^{0},t_{i}^{m}]$,
is to derive the optimal acceleration/deceleration which minimizes a convex
combination of its travel time and energy consumption. For each CAV $i$, we
define its \emph{information set} $Y_{i}(t)$, $t\in\lbrack t_{i}^{0},t_{i}%
^{f}]$, as
\[
Y_{i}(t)\triangleq\Big\{p_{i}(t),v_{i}(t),d_{i},s_{i}(t),\mathcal{I}%
_{i}\Big\},
\]
where $p_{i}(t),v_{i}(t)$ are the traveling distance and speed of CAV $i$
and $d_{i}$ indicates whether $i$ is making left
or right turn or going straight at the MZ. The fourth element in $Y_{i}(t)$ is
$s_{i}(t)=p_{k}(t)-p_{i}(t)$, the inter-vehicle distance between CAV $i$ and CAV $k$ which
is physically ahead of $i$ in the same lane (the index $k$ is made available
to $i$ by the coordinator). Based on this information the coordinator can also
evaluate%
\[
\mathcal{I}_{i}=\{e,s,l,o\}
\]
consisting of CAV indices as previously defined for the corresponding sets
$\{\mathcal{E}_{i}(t),$ $\mathcal{S}_{i}(t),$ $\mathcal{L}_{i}(t),$
$\mathcal{O}_{i}(t)\}$ known to the coordinator. The values of $e$, $s$, $l$,
$o$ are computed when CAV $i$ arrives at the CZ which is why $\mathcal{I}_{i}$
is treated as a time-invariant set. Since the coordinator is not involved in
any decision making process regarding vehicle control, we can formulate a
tractable decentralized problem, that can be solved on line by each CAV, as
follows:
\begin{equation}
\begin{aligned} &\min_{u_{i}\in U_{i}}\int_{t_{i}^{0}}^{t_{i}^{m}}[\gamma_{1}+\gamma_{2}u_{i}^{2}(t)]~dt \\ \text{s.t.}:&~\eqref{eq:model2}, (\ref{speed_accel constraints}),(\ref{rearend}),(\ref{Theorem1}), (\ref{tm_up_constraint}), \\ & ~p_i(t_i^0)=0, p_{i}(t_{i}^{m})=L \\ & \text{ and given }t_i^0, v_i(t_i^0), \label{eq:decentral} \end{aligned}
\end{equation}
where $\gamma_{1}=\beta$, $\gamma_{2}=\frac{(1-\beta)}%
{\bar{u}^{2}}$ with $\beta\in\lbrack0,1]$ a weight associated with the
importance of travel time relative to energy, and $\bar{u}=\max\{u_{i,\max
},\left\vert u_{i,\min}\right\vert \}$. In this manner, we have $\frac{u_{i}^{2}(t)}{\bar{u}^{2}%
}\in\lbrack0,1]$ and the overall objective is a well-defined convex
combination of the normalized travel time and normalized energy. Note that the
objective function in \eqref{eq:decentral} can be rewritten as%
\[
\gamma_{1}(t_{i}^{m}-t_{i}^{0})+\gamma_{2}\int_{t_{i}^{0}}^{t_{i}^{m}}%
u_{i}^{2}(t)dt
\]

\begin{customrem}
Unlike the problem considered in \cite{ZhangMalikopoulosCassandras2016} where
$t_{i}^{m}$ was obtained a priori to optimize travel times, here the optimal
travel time is determined by the solution of the problem itself. Similarly,
the terminal speed $v_{i}^{m}$ is also obtained from the optimal control
problem. An alternative formulation is to pre-specify a \textquotedblleft
target\textquotedblright\ $v_{i}^{m}$ or to include a penalty term on the
deviation of the actual terminal speed from $v_{i}^{m}$.
\end{customrem}

\subsection{Problem Decomposition \label{CZ_problem_decomposition}}

In order to simplify notation, we rewrite \eqref{eq:decentral} as
\begin{equation}
\begin{aligned} & \int_{t_{i}^{0}}^{t_{i}^{m}}\left[\gamma+ \frac{1}{2}u_{i}^{2}(t)\right] ~dt\\ \text{s.t.}: & ~\eqref{eq:model2}, (\ref{speed_accel constraints}),(\ref{rearend}),(\ref{Theorem1}), (\ref{tm_up_constraint}), \\ & ~p_i(t_i^0)=0, p_{i}(t_{i}^{m})=L, \\ & \text{ and given }t_i^0, v_i(t_i^0), \label{eq:p_ori} \end{aligned}
\end{equation}
where $\gamma=\frac{\gamma_{1}}{2\gamma_{2}}$. Note that \eqref{Theorem1} and
\eqref{tm_up_constraint} are constraints applied to the terminal time
$t_{i}^{m}$. In order to efficiently obtain analytical solutions on line, we
proceed with a two-step approach. The first step is to consider the following
problem by relaxing the terminal time constraints in \eqref{eq:p_ori}, that
is,
\begin{equation}
\begin{aligned} P_0:~~& \int_{t_{i}^{0}}^{t_{i}^{m}}\left[\gamma+ \frac{1}{2}u_{i}^{2}(t)\right] ~dt\\ \text{s.t.}: &~ \eqref{eq:model2}, (\ref{speed_accel constraints}),(\ref{rearend}), \\ & ~p_i(t_i^0)=0, p_{i}(t_{i}^{m})=L, \\ & \text{ and given }t_i^0, v_i(t_i^0). \label{P_0} \end{aligned}
\end{equation}
We denote the problem above as $P_{0}$. After solving $P_{0}$, the terminal
time obtained is denoted as $t_{i}^{m_{0}}$. The second step is to check
$t_{i}^{m_{0}}$ and $(i)$ If $t_{i}^{m_{0}}$ satisfies both \eqref{Theorem1}
and \eqref{tm_up_constraint}, then $t_{i}^{m_{0}}$ is the optimal terminal
time, $(ii)$ If $t_{i}^{m_{0}}$ violates \eqref{Theorem1}, then solve problem
$P_{1}$ by adding to $P_{0}$ the terminal time lower bound constraint from
Theorem 1, i.e., $t_{i}^{m}=\max\{t_{i}^{L},t_{e}^{f}+\frac{\delta}{v_{e}^{f}%
}-\Delta_{i},t_{s}^{m}+\Delta_{s}^{\delta},t_{s}^{f}-\Delta_{i},t_{l}%
^{f},t_{o}^{f}-\Delta_{i}\}$, $(iii)$ If $t_{i}^{m_{0}}$ violates
\eqref{tm_up_constraint} , then solve problem $P_{2}$ by adding to $P_{0}$ the
upper bound terminal time constraint (\ref{tm_up_constraint}) $t_{i}^{m}%
=t_{i}^{U}$. Note that if $t_{i}^{L}>t_{i}^{U}$, i.e., the lower bound on
$t_{i}^{m}$ is higher than its upper bound, the problem is obviously infeasible.

\subsection{Analytical Solution \label{CZ_solution}}

Given the objective function of problem $P_{0}$, the Hamiltonian is
\begin{equation}
\begin{aligned} H_i(p_i, v_i, u_i, \lambda_i, t) = \gamma + \frac{1}{2} u_i^2(t) + \lambda_i^pv_i(t) + \lambda_i^vu_i(t) \end{aligned}\label{hamil}%
\end{equation}
and the Lagrangian with constraints directly adjoined is
\begin{equation}
\begin{aligned} L_i(p_i, v_i, u_i, \lambda_i, \mu_i, \nu_i, t) &= H_i(p_i, v_i, u_i, \lambda_i, t) + \mu_i g_i(u_i, t) \\ & + \nu_i h_i(p_i, v_i, t) \end{aligned}\label{lag}%
\end{equation}
where (omitting time arguments for simplicity) $\lambda_{i}=[\lambda_{i}%
^{p},\lambda_{i}^{v}]^{T}\in\mathbb{R}^{2}$ is the costate vector,
$g_{i}(u_{i},t)\leq0$ and $h_{i}(p_{i},v_{i},t)\leq0$ represents the control
constraints and state constraints respectively, and $\mu_{i}=[\mu_{i}^{a}%
,\mu_{i}^{b}]^{T}\in\mathbb{R}^{2}$, $\nu_{i}=[\nu_{i}^{c},\nu_{i}^{d},\nu
_{i}^{s}]^{T}\in\mathbb{R}^{3}$ are Lagrange multipliers with
\begin{equation}
\mu_{i}^{a}=\left\{
\begin{array}
[c]{ll}%
>0, & \mbox{$u_{i}(t) - u_{max} =0$},\\
=0, & \mbox{$u_{i}(t) - u_{max} <0$},
\end{array}
\right.  \label{eq:17a}%
\end{equation}%
\begin{equation}
\mu_{i}^{b}=\left\{
\begin{array}
[c]{ll}%
>0, & \mbox{$u_{min} - u_{i}(t) =0$},\\
=0, & \mbox{$u_{min} - u_{i}(t)<0$},
\end{array}
\right.  \label{eq:17b}%
\end{equation}%
\begin{equation}
\nu_{i}^{c}=\left\{
\begin{array}
[c]{ll}%
>0, & \mbox{$v_{i}(t) - v_{max} =0$},\\
=0, & \mbox{$v_{i}(t) - v_{max}<0$},
\end{array}
\right.  \label{eq:17c}%
\end{equation}%
\begin{equation}
\nu_{i}^{d}=\left\{
\begin{array}
[c]{ll}%
>0, & \mbox{$v_{min} - v_{i}(t)=0$},\\
=0, & \mbox{$v_{min} - v_{i}(t)<0$}.
\end{array}
\right.  \label{eq:17d}%
\end{equation}%
\begin{equation}
\nu_{i}^{s}=\left\{
\begin{array}
[c]{ll}%
>0, & \mbox{$p_i(t) + \delta - p_k(t) = 0$},\\
=0, & \mbox{$p_i(t) + \delta - p_k(t) <0$}.
\end{array}
\right.  \label{eq:17e}%
\end{equation}
where $k$ is the CAV physically ahead of $i$ in the same lane and its position
$p_{k}(t)$ is known to $i$ through the coordinator (or through on-board sensors).

The Euler-Lagrange equations become
\begin{equation}
\dot{\lambda}_{i}^{p}(t)=-\frac{\partial L_{i}}{\partial p_{i}}=\left\{
\begin{array}
[c]{ll}%
0, & \mbox{$p_i(t) + \delta - p_k(t) < 0$},\\
-\nu_{i}^{s}, & \mbox{$p_i(t) + \delta - p_k(t) = 0$},
\end{array}
\right.  \label{eq:lambda_p}%
\end{equation}
and
\begin{equation}
\dot{\lambda}_{i}^{v}(t)=-\frac{\partial L_{i}}{\partial v_{i}}=\left\{
\begin{array}
[c]{ll}%
-\lambda_{i}^{p}(t), & \mbox{$v_{i}(t) - v_{max} <0$}~\text{and}\\
& \mbox{$v_{min} - v_{i}(t)<0$},\\
-\lambda_{i}^{p}(t)-\nu_{i}^{c}, & \mbox{$v_{i}(t) - v_{max} =0$},\\
-\lambda_{i}^{p}(t)+\nu_{i}^{d}, & \mbox{$v_{min} - v_{i}(t)=0$}.
\end{array}
\right.  \label{eq:lambda_v}%
\end{equation}
Since the terminal time $t_{i}^{m}$ and terminal speed $v_{i}^{m}$ are not
fixed, we have the transversality conditions
\begin{equation}
\lambda_{i}^{v}(t_{i}^{m})=0,\text{ \ \ }H_{i}(t_{i}^{m})=0
\label{eq:costate_trans}%
\end{equation}
Note that if we need to solve $P_{1}$ and $P_{2}$, then the second equation in
(\ref{eq:costate_trans}) no longer holds; instead, $t_{i}^{m}$ is fixed to its
upper or lower bound as described above. In addition, there also exist the
state boundary conditions $p_{i}(t_{i}^{0})=0$, $p_{i}(t_{i}^{m})=L$,
$v_{i}(t_{i}^{0})=v_{i}^{0}$, given the initial time $t_{i}^{0}$ and the initial speed $v_i^0$.

The necessary condition for optimality is
\begin{equation}
\frac{\partial L_{i}}{\partial u_{i}}=\left\{
\begin{array}
[c]{ll}%
u_{i}(t)+\lambda_{i}^{v}(t), & \mbox{$u_{i}(t) - u_{max} <0$}~\text{and}\\
& \mbox{$u_{min} - u_{i}(t)<0$},\\
u_{i}(t)+\lambda_{i}^{v}(t)+\mu_{i}^{a}, & \mbox{$u_{i}(t) - u_{max} =0$},\\
u_{i}(t)+\lambda_{i}^{v}(t)-\mu_{i}^{b}, & \mbox{$u_{min} - u_{i}(t)=0$}.
\end{array}
\right.  \label{eq:optimum}%
\end{equation}
A complete solution of this problem requires that constrained and
unconstrained arcs of an optimal trajectory are pieced together to satisfy all
conditions (\ref{eq:17a}) through (\ref{eq:optimum}). This includes the five
constraints (three pure-state constraints, two control constraints) in
(\ref{eq:17a}) through (\ref{eq:17e}). While there are many different cases
that can occur, the nature of the optimal solution rules out the possibility
of several cases. In what follows, we provide a complete analysis of the case
where no constraints are active and of the case where the safety constraint
$p_{i}(t)+\delta-p_{k}(t)\leq0$ is the only active one. A discussion of the remaining cases can be found in Appendix \ref{appendix:CZ}.


\subsection{Unconstrained Optimal Control Analysis \label{no_active}}

For problem $P_{0}$, the terminal time is free whereas for $P_{1}$ and $P_{2}$
the terminal time is fixed. Thus, we provide the analysis for each of these
two cases.

\subsubsection{Free Terminal Time}

\label{no_active_free} When the state and control constraints are inactive, we
have $\mu_{i}^{a}=\mu_{i}^{b}=\nu_{i}^{c}=\nu_{i}^{d}=\nu_{i}^{s}=0$. The
Lagrangian \eqref{lag} becomes $L_{i}(p,v,u,\lambda,\mu,\nu,t)=H_{i}%
(p,v,u,\lambda,t)$ and \eqref{eq:optimum} reduces to $\frac{\partial L_{i}%
}{\partial u_{i}}=u_{i}(t)+\lambda_{i}^{v}=0$, which leads to
\begin{equation}
u_{i}(t)=-\lambda_{i}^{v}(t).\label{optimum_1}%
\end{equation}
Since $\nu_{i}^{s}=0$, \eqref{eq:lambda_p} becomes $\dot{\lambda}_{i}%
^{p}(t)=-\frac{\partial L_{i}}{\partial p_{i}}=0$ which leads to
\begin{equation}
\lambda_{i}^{p}=a_{i}\label{lambdaip}%
\end{equation}
where $a_{i}$ is a constant. Since $\nu_{i}^{c}=\nu_{i}^{d}=0$,
\eqref{eq:lambda_v} becomes $\dot{\lambda}_{i}^{v}(t)=-\frac{\partial L_{i}%
}{\partial v_{i}}=-\lambda_{i}^{p}$. Since $\lambda_{i}^{p}=a_{i}$, we have
\begin{equation}
\lambda_{i}^{v}(t)=-a_{i}t-b_{i}\label{lambdaiv}%
\end{equation}
where $b_{i}$ is a constant. We can now obtain a complete analytical soloution
of $P_{0}$ as follows.

\begin{customthm}
The optimal trajectory for problem $P_{0}$ is given by
\begin{align}
u_{i}^{\ast}(t)  &  =a_{i}t+b_{i}\label{ui*}\\
v_{i}^{\ast}(t)  &  =\frac{1}{2}a_{i}t^{2}+b_{i}t+c_{i}\label{vi*}\\
p_{i}^{\ast}(t)  &  =\frac{1}{6}a_{i}t^{3}+\frac{1}{2}b_{i}t^{2}+c_{i}t+d_{i}
\label{pi*}%
\end{align}
for $t\in\lbrack t_{i}^{0},t_{i}^{m^{\ast}}]$ where $a_{i}$, $b_{i}$, $c_{i}$
and $d_{i}$ are constants determined along with $t_{i}^{m^{\ast}}$ through
\begin{subequations}
\begin{align}
\frac{1}{6}a_{i} \cdot (t_{i}^{0})^{3} + \frac{1}{2}b_{i} \cdot  (t_{i}^{0})^{2} + c_{i}
t_{i}^{0} + d_{i}  &  = 0
\label{P0_solution:p0}\\
\frac{1}{2}a_{i} \cdot  (t_{i}^{0})^{2} + b_{i} t_{i}^{0} + c_{i}  &  = v_{i}%
^{0}\label{P0_solution:v0}\\
\frac{1}{6}a_{i} \cdot  (t_{i}^{m})^{3} + \frac{1}{2}b_{i} \cdot  (t_{i}^{m})^{2} + c_{i}
t_{i}^{m} + d_{i}  &  = L\label{P0_solution:pm}\\
a_{i} t_{i}^{m} + b_{i}  & = 0\label{P0_solution:vm_free}\\
\gamma- \frac{1}{2}b_{i}^{2} + a_{i} c_{i}  &  = 0\label{P0_solution:tran}%
\end{align}
\label{P0_solution}
\end{subequations}
\end{customthm}

\emph{Proof. } The optimal control in (\ref{ui*}) follows from
\eqref{optimum_1} and (\ref{lambdaiv}). Using \eqref{ui*} in the system
dynamics \eqref{eq:model2}, we then derive (\ref{vi*}) and (\ref{pi*}). Next,
(\ref{P0_solution:p0}) through (\ref{P0_solution:pm}) follow from the boundary
conditions $p_{i}(t_{i}^{0})=0$, $v_{i}(t_{i}^{0})=v_{i}^{0}$, $p_{i}%
(t_{i}^{m})=L$ and (\ref{P0_solution:vm_free}) follows from $\lambda_{i}%
^{v}(t_{i}^{m})=0$ in \eqref{eq:costate_trans} and from (\ref{ui*}). The last
equation follows from $H_{i}(t_{i}^{m})=0$ in \eqref{eq:costate_trans}:
\begin{align*}
&  \gamma+\frac{1}{2}(u_{i}^{\ast}(t_{i}^{m}))^{2}+a_{i}v_{i}^{\ast}(t_{i}%
^{m})-(u_{i}^{\ast}(t_{i}^{m}))^{2}\\
&  =\gamma-\frac{1}{2}(a_{i}t_{i}^{m}+b_{i})^{2}+a_{i}(\frac{1}{2}a_{i}%
(t_{i}^{m})^{2}+b_{i}t_{i}^{m}+c_{i})\\
&  =\gamma-\frac{1}{2}b_{i}^{2}+a_{i}c_{i}=0
\end{align*}
using (\ref{hamil}), (\ref{lambdaip}), (\ref{lambdaiv}), (\ref{ui*}) and (\ref{vi*}). \hfill$\blacksquare$

Thus, a complete solution of $P_{0}$ boils down to solving the five equations
in (\ref{P0_solution}). A typical simulation example of this case can be found
in Section \ref{analytic}.

The next two results establish a basic property of the optimal control, i.e.,
it is non-negative and non-increasing, and the fact that two of the
constraints in (\ref{speed_accel constraints}) cannot be active.


\begin{customlem}
For the unconstrained problem with free terminal time, the optimal control is
non-negative, i.e., $u_{i}^{\ast}(t)\geq0$, and monotonically non-increasing
\label{lem:u_nonnegative}
\end{customlem}

\emph{Proof. } To prove that $u_{i}^{\ast}(t)\geq0$, let us assume that the
optimal control includes an interval that is negative, i.e., $u_{i}^{\ast
}(t)<0$ for $t\in\lbrack t_{1},t_{2}]$, and $u_{i}^{\ast}(t)\geq0$ for
$t\in\lbrack t_{i}^{0},t_{1}]\cup\lbrack t_{2},t_{i}^{m}]$. Next, we construct
another control $u_{i}^{c}(t)$ which is the same as $u_{i}^{\ast}(t)$ except
that $u_{i}^{c}(t)=0$ for $t\in\lbrack t_{1},t_{2}]$. Let us first prove that
$u_{i}^{c}(t)$ is feasible. Note that $u_{i}^{c}(t)=u_{i}^{\ast}(t)$ for
$t\in\lbrack t_{i}^{0},t_{1}]\cup\lbrack t_{2},t_{i}^{m}]$, hence, $u_{i}%
^{c}(t)$ must be feasible for $t\in\lbrack t_{i}^{0},t_{1}]$. Moreover, for
$t\in\lbrack t_{1},t_{2}]$, $u_{i}^{c}(t)=0$ does not violate any control
constraints. In addition, since we have $v_{min}<v_{i}^{c}(t_{1})<v_{max}$ and
$u_{i}^{c}(t)=0$ for $t\in\lbrack t_{1},t_{2}]$, it follows that $v_{i}%
^{c}(t)=v_{i}(t_{1})$ and $v_{min}<v_{i}^{c}(t)<v_{max}$ for $t\in\lbrack
t_{1},t_{2}]$. Hence, the speed constraints are satisfied as well. By
assumption, the safety constraint (\ref{rearend}) is inactive over $[t_{i}%
^{0},t_{1}]$. If, under $u_{i}^{c}(t)$, the safety constraint remains inactive
over $[t_{1},t_{2}]$, then $u_{i}^{c}(t)$ is feasible. Otherwise, since
(\ref{rearend}) is inactive under $u_{i}^{\ast}(t)$, there exists some
$\epsilon>0$ such that it remains inactive over $[t_{1},t_{1}+\epsilon]$, in
which case we modify $u_{i}^{c}(t)$ so that $u_{i}^{c}(t)=0$ for $t\in\lbrack
t_{1},t_{1}+\epsilon]$ and $u_{i}^{c}(t)=u_{i}^{\ast}(t)$ for $t\in
(t_{1}+\epsilon,t_{2}]$ which ensures $u_{i}^{c}(t)$ is feasible over
$[t_{1},t_{2}]$. Next, since $u_{i}^{\ast}(t)<0$ and $u_{i}^{c}(t)=0$ for
$t\in\lbrack t_{1},t_{1}+\epsilon]$, we have $v_{i}^{\ast}(t)<v_{i}^{c}(t)$
for $t\in\lbrack t_{1},t_{2}]$ and for $t\in\lbrack t_{2},t_{i}^{m}]$. This
implies that it is possible that $v_{i}^{c}(t_{3})=v_{max}$ for some $t_{3}%
\in(t_{2},t_{i}^{m}]$. If that happens, then we modify $u_{i}^{c}(t)$ so that
$u_{i}^{c}(t)=0$, hence $v_{i}^{c}(t)=v_{max}$, for $t\in\lbrack t_{3}%
,t_{i}^{m}]$ and $u_{i}^{c}(t)$ is feasible over all $t\in\lbrack t_{i}%
^{0},t_{i}^{m}]$. Since $v_{i}^{\ast}(t)<v_{i}^{c}(t)$ for $t\in\lbrack
t_{1},t_{2}]$ and $v_{i}^{\ast}(t)\leq v_{i}^{c}(t)$ for $t\in\lbrack
t_{i}^{0},t_{i}^{m}]$, it follows that $t_{i}^{m^{\ast}}>t_{i}^{m^{c}}$, or
$t_{i}^{m^{\ast}}-t_{i}^{0}>t_{i}^{m^{c}}-t_{i}^{0}$, which indicates that the
optimal control $u_{i}^{\ast}(t)$ leads to a longer travel time compared to
$u_{i}^{c}(t)$. Denoting the energy consumption in \eqref{eq:p_ori} as
$J_{u}(u_{i}(t))$ and allowing for the possibility that there exists $t_{3}$
such that $t_{2}<t_{3}\leq t_{i}^{m}$ (otherwise, $t_{3}=t_{i}^{m}$), we have%
\begin{align*}
J_{u}(u_{i}^{\ast}(t)) &  =\int_{t_{i}^{0}}^{t_{1}}\frac{1}{2}(u_{i}^{\ast
})^{2}dt+\int_{t_{1}}^{t_{2}}\frac{1}{2}(u_{i}^{\ast})^{2}dt\\
&  +\int_{t_{2}}^{t_{3}}\frac{1}{2}(u_{i}^{\ast})^{2}dt+\int_{t_{3}}%
^{t_{i}^{m}}\frac{1}{2}(u_{i}^{\ast})^{2}dt\\
&  >\int_{t_{i}^{0}}^{t_{1}}\frac{1}{2}(u_{i}^{\ast})^{2}dt+\int_{t_{2}%
}^{t_{3}}\frac{1}{2}(u_{i}^{\ast})^{2}dt=J_{u}(u_{i}^{c}(t))
\end{align*}
where we have used the fact that $u_{i}^{c}(t)=0$ for at least $t\in\lbrack
t_{1},t_{1}+\epsilon]$ and $u_{i}^{c}(t)=u_{i}^{\ast}(t)$ for $t\in\lbrack
t_{i}^{0},t_{1}]\cup\lbrack t_{2},t_{3}]$. In view of the objective function
\eqref{eq:decentral}, since we have shown that $t_{i}^{m^{\ast}}-t_{i}%
^{0}>t_{i}^{c}-t_{i}^{0}$ and $J_{u}(u_{i}^{\ast}(t))>J_{u}(u_{i}%
^{c}(t))$, we can see that the optimal control $u_{i}^{\ast}(t)$ incurs a
larger cost compared to $u_{i}^{c}(t)$, which contradicts the assumption that
$u_{i}^{\ast}(t)$ is the optimal control. Therefore, $u_{i}^{\ast}(t)$ cannot
be negative at any $t\in\lbrack t_{i}^{0},t_{i}^{m}]$ and we have $u_{i}%
^{\ast}(t)\geq0$ as long as the optimal trajectory is unconstrained.

To prove that $u_{i}^{\ast}(t)$ is monotonically non-increasing, observe that
$u_{i}^{\ast}(t_{i}^{m})=0$ from (\ref{P0_solution:vm_free}). Since the
optimal control is a linear function of time as indicated by \eqref{ui*} and
non-negative as established above, it follows that $u_{i}^{\ast}(t)$ is
monotonically non-increasing. \hfill$\blacksquare$

\begin{customlem}
For the unconstrained problem with free terminal time, it is not possible for
constraints $v_{min}- v_{i}(t) \leq0$ and/or $u_{min} - u_{i}(t) \leq0$ to
become active. \label{lem:rule_out_dec}
\end{customlem}

\emph{Proof. } By Lemma \ref{lem:u_nonnegative}, $u_{i}^{\ast}(t)\geq0$.
Hence, $u_{min}-u_{i}(t)\leq0$ cannot become active. Since, by assumption,
$v_{min}<v_{i}^{0}<v_{max}$, $v_{i}^{\ast}(t)$ cannot reach $v_{min}$ and the
constraint $v_{min}-v_{i}(t)\leq0$ cannot become active. \hfill$\blacksquare$


\subsubsection{Fixed Terminal Time}

\label{no_active_fixed_tm} If the terminal time $t_{i}^{m^{\ast}}$ obtained
from solving $P_{0}$ turns out to violate \eqref{Theorem1} or
\eqref{tm_up_constraint}, then, as described in the two-step approach earlier,
we need to solve either $P_{1}$ or $P_{2}$ by setting $t_{i}^{m}$ to a fixed
value which is either the lower bound in \eqref{Theorem1} or the upper bound
in \eqref{tm_up_constraint}. Therefore, the transversality condition
$H_{i}(t_{i}^{m})=0$, i.e., the fifth equation in (\ref{P0_solution}), no
longer holds and the solution of this problem reduces to
\begin{equation}
\left[
\begin{array}
[c]{cccc}%
\frac{1}{6}(t_{i}^{0})^{3} & \frac{1}{2}(t_{i}^{0})^{2} & t_{i}^{0} & 1\\
\frac{1}{2}(t_{i}^{0})^{2} & t_{i}^{0} & 1 & 0\\
\frac{1}{6}(t_{i}^{m})^{3} & \frac{1}{2}(t_{i}^{m})^{2} & t_{i}^{m} & 1\\
t_{i}^{m} & 1 & 0 & 0
\end{array}
\right]  .\left[
\begin{array}
[c]{c}%
a_{i}\\
b_{i}\\
c_{i}\\
d_{i}%
\end{array}
\right]  =\left[
\begin{array}
[c]{c}%
0\\
v_{i}^{0}\\
L\\
0
\end{array}
\right]  \label{fixed_term_time_solution}%
\end{equation}
which yields the four parameters $a_{i}$, $b_{i}$, $c_{i}$, $d_{i}$ from a
simple system of linear equations. A typical simulation example of this case when \eqref{Theorem1} is violated can be found in Section \ref{analytic}.

With the terminal time fixed, Lemma \ref{lem:u_nonnegative} needs to be
modified as follows.

\begin{customlem}
For the unconstrained problem with fixed terminal time, the optimal control
must be either monotonically non-increasing and $u_{i}^{\ast}(t)\geq0$, or
monotonically non-decreasing and $u_{i}^{\ast}(t)\leq0.$ \label{lem:u_fixed}
\end{customlem}

\emph{Proof. } Due to the linearity of the unconstrained optimal control in
\eqref{ui*}, $u_{i}^{\ast}(t)$ must be either non-increasing or non-decreasing
over all $t\in\lbrack t_{i}^{0},t_{i}^{m}]$. With the terminal time fixed, the
last equation in (\ref{fixed_term_time_solution}) gives $u_{i}^{\ast}%
(t_{i}^{m})=0$. Therefore, when $u_{i}^{\ast}(t)$ is monotonically
non-increasing, we have $u_{i}^{\ast}(t)\geq0$; when $u_{i}^{\ast}(t)$ is
monotonically non-decreasing, we have $u_{i}^{\ast}(t)\leq0$. \hfill
$\blacksquare$

\subsection{Constrained Optimal Control Analysis \label{active}}

Checking whether the optimal solution of the unconstrained problem $P_{0}$,
$P_{1}$ or $P_{2}$ violates any of the constraints (\ref{eq:17a}) through
(\ref{eq:17e}) is easily accomplished since the unconstrained optimal control
in \eqref{ui*} is a linear function of time and the optimal speed is a
quadratic function of time. When this happens, we must check whether there
exists a nonempty feasible control set. One approach followed in earlier work
\cite{Zhang2016} is to identify the set of all initial conditions $(t_{i}%
^{0},v_{i}^{0})$ such that no constraint is violated over $[t_{i}^{0}%
,t_{i}^{m}]$ or at least some of the constraints are not violated while the
rest are explicitly dealt with through the Lagrangian in (\ref{lag}). As shown
in \cite{Zhang2016}, it is possible to define a Feasibility Enforcement Zone
(FEZ) which precedes the CZ such that each CAV is controlled over the FEZ so
as to reach a feasible initial condition when reaching the CZ. Here, however,
we proceed differently by following a direct approach through which we derive
explicit solutions for any feasible optimal constrained trajectory. In so
doing, we can also explicitly identify when an optimal solution is infeasible
under initial conditions $(t_{i}^{0},v_{i}^{0})$.

When the optimal solution of the unconstrained problem violates a constraint,
we need to re-solve the problem by identifying an optimal trajectory that
includes unconstrained arcs pieced together with one or more constrained arcs
such that all necessary conditions for optimality are satisfied. For a control
constraint of the form $g_{i}(u_{i},t)\leq0$ as in (\ref{eq:17a}%
)-(\ref{eq:17b}), the optimal control on a constrained arc can be simply
obtained by solving $g_{i}(u_{i},t)=0$. The remaining constraints
(\ref{eq:17c})-(\ref{eq:17e}) in our problem are pure state constraints of the
form $h_{i}(x_{i},t)\leq0$. In this case (see \cite{bryson1975optimal}), we
define the tangency constraints
\begin{equation}
N_{i}(x_{i},t)\triangleq\lbrack h_{i}(x_{i},t)\text{ }h_{i}^{(1)}%
(x_{i},t)\text{ }\cdots\text{ }h_{i}^{(q-1)}(x_{i},t)]^{T}=0,
\label{tangent_vp}%
\end{equation}
where $h_{i}^{(k)}(x_{i},t)$ is the $k$th time derivative and $q$ derivatives
are taken until we obtain an expression that explicitly depends on the control
$u_{i}$ so that
\begin{equation}
h_{i}^{(q)}(x_{i},t)=0. \label{tangent_u}%
\end{equation}
At the junction points of constrained and unconstrained arcs, the costate and
Hamiltonian trajectories may have discontinuities. This can be determined
using the following jump conditions \cite{bryson1975optimal}, where $\tau$
denotes a junction point and $\tau^{-}$,$\tau^{+}$ denote the left-hand side
and the right-hand side limits, respectively:
\begin{equation}
\begin{aligned} \lambda_i(\tau^-) &= \lambda_i(\tau^+) + \pi_i^T \frac{\partial N_i(x_i,t)}{\partial x_i},\\ H_i(\tau^-) &= H_i(\tau^+) - \pi_i^T \frac{\partial N_i(x_i,t)}{\partial t}. \end{aligned} \label{ipm_disc}%
\end{equation}
where $N_{i}(x_{i},t)$ is the $q$-dimensional vector in (\ref{tangent_vp}) and
$\pi_{i}$ is a $q$-dimensional vector of constant Lagrange multipliers
satisfying $\pi_{i}^{T}N_{i}(x_{i},t)=0$ and $\pi_{i}\geq0$, $i=1,\ldots,q$.
Consequently, the optimal control $u_{i}^{\ast}(t)$ may or may not be
continuous at junction points.

In what follows, we concentrate on the safety constraint (\ref{rearend}) which
is the most challenging to deal with. In this case, we have $\mu_{i}^{a}%
=\mu_{i}^{b}=\nu_{i}^{c}=\nu_{i}^{d}=0$. The remaining constraints are
discussed in the Appendix. Thus, we set $h_{i}(p_{i},t)=p_{i}+\delta
-p_{k}^{\ast}(t)$ where we observe that $p_{k}^{\ast}(t)$ is a known explicit
function of time given by the optimal position trajectory of CAV $k$ specified
in (\ref{P0_solution}) or (\ref{fixed_term_time_solution}) since, upon arrival
of CAV $i$ at the CZ, the optimal solution of the problem associated with
$k<i$ has already been fully determined. Moreover, $h_{i}^{(1)}(p_{i}%
,t)=v_{i}-\frac{\partial p_{k}^{\ast}(t)}{\partial t}=v_{i}-v_{k}^{\ast}(t)$
where $v_{k}^{\ast}(t)$ is also an explicit function of time in
(\ref{P0_solution}) or (\ref{fixed_term_time_solution}) and $h_{i}^{(2)}%
(p_{i},t)=u_{i}-\frac{\partial v_{k}^{\ast}(t)}{\partial t}=u_{i}-u_{k}^{\ast
}(t)$, hence the optimal control on the constrained arc is given by
$u_{i}^{\ast}(t)=u_{k}^{\ast}(t)$.
%

The following result establishes the continuity property of the optimal
control when the trajectory enters a constrained arc where $p_{i}%
(t)+\delta-p_{k}^{\ast}(t)=0$.


\begin{customthm}
The optimal control $u_{i}^{\ast}(t)$ is continuous at the junction $\tau$ of
the unconstrained and safety-constrained arcs, i.e., $u_{i}^{\ast}(\tau^{-}%
)=u_{i}^{\ast}(\tau^{+})$. \label{safety:u_cont}
\end{customthm}


\emph{Proof. } By assumption, the rear-end safety constraint is not active at
$t_{i}^{0}$. Hence, when the safety constraint $p_{i}(t)+\delta-p_{k}^{\ast
}(t)\leq0$ becomes active, $\tau$ is the entry time of the constrained arc,
and the jump conditions in \eqref{ipm_disc} become%
\begin{align*}
\lambda_{i}^{p}(\tau^{-}) &  =\lambda_{i}^{p}(\tau^{+})+\pi_{i}^{p}%
\frac{\partial}{\partial p_{i}}[p_{i}+\delta-p_{k}^{\ast}(t)]\\
\lambda_{i}^{v}(\tau^{-}) &  =\lambda_{i}^{v}(\tau^{+})+\pi_{i}^{v}%
\frac{\partial}{\partial v_{i}}[v_{i}-v_{k}^{\ast}(t)]\\
H_{i}(\tau^{-}) &  =H_{i}(\tau^{+})-\pi_{i}^{p}\frac{\partial}{\partial
t}[p_{i}+\delta-p_{k}^{\ast}(t)]-\pi_{i}^{v}\frac{\partial}{\partial t}%
[v_{i}-v_{k}^{\ast}(t)]
\end{align*}
where $\frac{\partial p_{k}^{\ast}(t)}{\partial t}=v_{k}^{\ast}(t)$ and
$\frac{\partial v_{k}^{\ast}(t)}{\partial t}=u_{k}^{\ast}(t)$ are explicit
functions of $t$ specified through (\ref{P0_solution}) or
(\ref{fixed_term_time_solution}). We assume that $u_{k}^{\ast}(t)$, $k<i$, is
continuous in $t$ so that, if we can establish that $u_{k}^{\ast}(t)$ is
continuous, then a simple iterative argument completes the proof. The
equations above become%
\begin{align*}
\lambda_{i}^{p}(\tau^{-}) &  =\lambda_{i}^{p}(\tau^{+})+\pi_{i}^{p},\text{
\ \ }\lambda_{i}^{v}(\tau^{-})=\lambda_{i}^{v}(\tau^{+})+\pi_{i}^{v},\\
H_{i}(\tau^{-}) &  =H_{i}(\tau^{+})+\pi_{i}^{p}v_{k}^{\ast}(t)+\pi_{i}%
^{v}u_{k}^{\ast}(t)
\end{align*}
For $t\geq\tau^{+}$, the tangency conditions (\ref{tangent_vp}%
)-(\ref{tangent_u}) with $q=2$ hold:%
\begin{align*}
p_{i}(t)+\delta-p_{k}^{\ast}(t) &  =0\\
v_{i}(t)-v_{k}^{\ast}(t) &  =0\\
u_{i}(t)-u_{k}^{\ast}(t) &  =0
\end{align*}
In addition, note that the position $p_{i}(t)$ and speed $v_{i}(t)$ are
continuous functions of $t$. Combining the equations above and recalling from
(\ref{hamil}) that $H_{i}(t)=\gamma+\frac{1}{2}u_{i}^{2}(t)+\lambda_{i}%
^{p}(t)v_{i}(t)+\lambda_{i}^{v}(t)u_{i}(t)$, we get%
\begin{align*}
&  \gamma+\frac{1}{2}u_{i}^{2}(\tau^{-})+\lambda_{i}^{p}(\tau^{-})v_{i}%
(\tau)+\lambda_{i}^{v}(\tau^{-})u_{i}(\tau^{-})\\
&  =\gamma+\frac{1}{2}u_{i}^{2}(\tau^{+})+\lambda_{i}^{p}(\tau^{+})v_{i}%
(\tau)+\lambda_{i}^{v}(\tau^{+})u_{i}(\tau^{+})\\
&  +\pi_{i}^{p}v_{k}^{\ast}(\tau)+\pi_{i}^{v}u_{k}^{\ast}(\tau)
\end{align*}
and since $u_{i}(\tau^{+})=u_{k}^{\ast}(\tau)$, it follows that%
\begin{gather*}
\frac{1}{2}u_{i}^{2}(\tau^{-})-\frac{1}{2}u_{i}^{2}(\tau^{+})+[\lambda_{i}%
^{p}(\tau^{-})-\lambda_{i}^{p}(\tau^{+})]v_{k}^{\ast}(\tau)\\
-\pi_{i}^{p}v_{k}^{\ast}(\tau)+\lambda_{i}^{v}(\tau^{-})u_{i}(\tau
^{-})-\lambda_{i}^{v}(\tau^{+})u_{i}(\tau^{+})-\pi_{i}^{v}u_{i}(\tau^{+})=0
\end{gather*}
which reduces to%
\begin{align*}
&  \frac{1}{2}u_{i}^{2}(\tau^{-})-\frac{1}{2}u_{i}^{2}(\tau^{+})+\lambda
_{i}^{v}(\tau^{-})[u_{i}(\tau^{-})-u_{i}(\tau^{+})]\\
&  =[u_{i}(\tau^{-})-u_{i}(\tau^{+})](\frac{1}{2}[u_{i}(\tau^{-})+u_{i}%
(\tau^{+})]+\lambda_{i}^{v}(\tau^{-}))=0
\end{align*}
Therefore, either $u_{i}(\tau^{-})-u_{i}(\tau^{+})=0$, or $\frac{1}{2}%
[u_{i}(\tau^{-})+u_{i}(\tau^{+})]+\lambda_{i}^{v}(\tau^{-})=0$. Assuming that
$u_{i}(\tau^{-})-u_{i}(\tau^{+})\neq0$, recall that at $\tau^{-}$ the
trajectory arc is unconstrained so that (\ref{optimum_1}) holds: $u_{i}%
(\tau^{-})=-\lambda_{i}^{v}(\tau^{-})$ and it follows that $u_{i}(\tau
^{-})-u_{i}(\tau^{+})=0$. We conclude that $u_{i}(t)$ is continuous at $\tau$
and the proof is complete. \hfill$\blacksquare$ \bigskip

Once an optimal trajectory for CAV $i$ enters the constrained arc
$p_{i}(t)+\delta-p_{k}^{\ast}(t)=0$, it may remain on this arc through the
terminal time $t_{i}^{m}$ or exit it at some point $\tau^{\prime}>\tau$ and
follow an unconstrained arc over $[\tau^{\prime},t_{i}^{m}]$. This depends on
whether such an exit point $\tau^{\prime}$ is feasible on an optimal
trajectory. More generally, it is possible that an optimal trajectory consists
of a sequence of alternating unconstrained and constrained arcs whose
feasibility needs to be checked. Thus, once we establish that an optimal
trajectory contains a constrained arc, there are two cases to consider.

\textbf{Case 1}: No exit point from the constrained arc. In this case, CAV $i$
remains on the constrained arc until it reaches the MZ and we have%
\begin{equation}
u_{i}^{\ast}(t)=\left\{
\begin{array}
[c]{cc}%
a_{i}t+b_{i} & t\in\lbrack t_{i}^{0},\tau]\\
u_{k}^{\ast}(t) & t\in(\tau,t_{i}^{m}]
\end{array}
\right.  \label{ui*_unconstr_arc}%
\end{equation}
Moreover, $v_{i}^{\ast}(t)$ is given by (\ref{vi*}) for $t\in\lbrack t_{i}%
^{0},\tau]$ and $v_{i}^{\ast}(t)=v_{k}^{\ast}(t)$ for $t\in(\tau,t_{i}^{m}]$;
$p_{i}^{\ast}(t)$ is given by (\ref{pi*}) for $t\in\lbrack t_{i}^{0},\tau]$
and $p_{i}^{\ast}(t)+\delta-p_{k}^{\ast}(t)=0$ for $t\in(\tau,t_{i}^{m}]$. The
constants $a_{i}$, $b_{i}$, $c_{i}$ and $d_{i}$ along with $\tau$ are
determined through
\begin{subequations}
\begin{align}
\frac{1}{6}a_{i}\cdot(t_{i}^{0})^{3}+\frac{1}{2}b_{i}\cdot(t_{i}^{0}%
)^{2}+c_{i}t_{i}^{0}+d_{i} &  =0\label{necessary_safety_active_2arcs:p0}\\
\frac{1}{2}a_{i}\cdot(t_{i}^{0})^{2}+b_{i}t_{i}^{0}+c_{i} &  =v_{i}%
^{0}\label{necessary_safety_active_2arcs:v0}\\
a_{i}\tau+b_{i} &  =u_{k}^{\ast}(\tau
)\label{necessary_safety_active_2arcs:u_tau}\\
\frac{1}{6}a_{i}\tau^{3}+\frac{1}{2}b_{i}\tau^{2}+c_{i}\tau+d_{i}+\delta
 &  =p_{k}^{\ast}(\tau)\label{necessary_safety_active_2arcs:p_tau}\\
\frac{1}{2}a_{i}\tau^{2}+b_{i}\tau+c_{i} &  =v_{k}^{\ast}(\tau
)\label{necessary_safety_active_2arcs:v_tau}%
\end{align}
\label{necessary_safety_active_2arcs}
\end{subequations}
The first two equations above are the same as in (\ref{P0_solution}) and
follow from the initial conditions., while
\eqref{necessary_safety_active_2arcs:u_tau} follows from
(\ref{ui*_unconstr_arc}) and Theorem \ref{safety:u_cont}. In addition,
\eqref{necessary_safety_active_2arcs:p_tau} follows from
\eqref{necessary_safety_active_2arcs:p0} when the safety constraint becomes
active, i.e., $p_{i}^{\ast}(\tau)+\delta-p_{k}^{\ast}(\tau)=0$, and
\eqref{necessary_safety_active_2arcs:v_tau} follows from $v_{i}^{\ast}%
(\tau)=v_{k}^{\ast}(\tau)$. Note that in this case the terminal
time $t_{i}^{m}$ is fixed and determined by CAV $s$ in \eqref{Theorem1}. A typical
simulation example of this case is given in Section \ref{analytic} (Fig.
\ref{safety:simulation_2arcs}).
\begin{customrem}
As noted in Section II, an alternative to the distance-based safety constraint
$p_{k}(t)-p_{i}(t)\geq\delta$ is the speed-based safety constraint
\cite{Rajamani2000} $p_{k}(t)-p_{i}(t)\geq\varphi v_{i}(t)+\delta$. While the
expression of the analytical solutions become more complicated, the approach
for deriving all necessary conditions is the same as described above.
\end{customrem}




\textbf{Case 2}: There exists an exit point from the constrained arc. In this
case, letting $\tau_{1}$ denote the entry point to the constrained arc and
$\tau_{2}$ the exit point, there are two subcases to consider: $(i)$ when the
terminal time $t_{i}^{m}$ is free, and $(ii)$ when the terminal time is fixed.
When the terminal time is free, the transversality condition
\eqref{eq:costate_trans} holds, and a solution is obtained through the system
of five equations in \eqref{necessary_safety_active_2arcs} with $\tau_{1}$
replacing $\tau$ along with the following equations:%
\begin{subequations}
\begin{align}
e_{i}\tau_{2}+r_{i} &  =u_{k}^{\ast}(\tau_{2}%
)\label{necessary_safety_active_3arcs_free:u_tau}\\
\frac{1}{6}e_{i}\tau_{2}^{3}+\frac{1}{2}r_{i}\tau_{2}^{2}+q_{i}\tau_{2}%
+m_{i}+\delta &
=p_{k}^{\ast}(\tau_{2})\label{necessary_safety_active_3arcs_free:p_tau}\\
\frac{1}{2}e_{i}\tau_{2}^{2}+r_{i}\tau_{2}+q_{i} &  =v_{k}^{\ast}(\tau
_{2})\label{necessary_safety_active_3arcs_free:v_tau}\\
\frac{1}{6}e_{i}\cdot(t_{i}^{m})^{3}+\frac{1}{2}r_{i}\cdot(t_{i}^{m}%
)^{2}+q_{i}t_{i}^{m}+m_{i} &  =L\label{necessary_safety_active_3arcs_free:pm}%
\\
e_{i}t_{i}^{m}+r_{i} &  =0\label{necessary_safety_active_3arcs_free:vm_free}\\
\gamma-\frac{1}{2}r_{i}^{2}+e_{i}q_{i} &
=0.\label{necessary_safety_active_3arcs_free:tm_free}%
\end{align}
\label{necessary_safety_active_3arcs_free}
\end{subequations}
The first equation above follows from the fact that $[\tau_{2},t_{i}^{m}]$ is
an unconstrained arc so that (\ref{ui*}) applies but with new constans $e_{i}%
$, $r_{i}$ and from $u_{i}^{\ast}(\tau_{2}^{-})=u_{i}^{\ast}(\tau_{2}%
^{+})=u_{k}^{\ast}(\tau_{2})$;
\eqref{necessary_safety_active_3arcs_free:p_tau} follows from the constraint
$p_{i}^{\ast}(\tau_{2})+\delta-p_{k}^{\ast}(\tau_{2})=0$ and
\eqref{necessary_safety_active_3arcs_free:v_tau} from $v_{i}^{\ast}(\tau
_{2})=v_{k}^{\ast}(\tau_{2})$. Next,
\eqref{necessary_safety_active_3arcs_free:pm} is the boundary constraint
$p_{i}^{\ast}(t_{i}^{m})=L$, while
\eqref{necessary_safety_active_3arcs_free:vm_free} and
\eqref{necessary_safety_active_3arcs_free:tm_free} are the transversality
conditions similar to the last two equations in (\ref{P0_solution}). A
simulation example of this case is given in Section \ref{analytic} (Fig.
\ref{safety:simulation_3arcs}). In the case where the terminal time $t_{i}%
^{m}$ is fixed (as in problems $P_{1}$ and $P_{2}$) we simply remove the
transversality condition in the last equation above.
\begin{customrem}
Note that $a_{i}$, $b_{i}$, $c_{i}$, $d_{i}$ and $\tau_{1}$ can be determined
from the five equations in \eqref{necessary_safety_active_2arcs} independently
from $e_{i}$, $r_{i}$, $q_{i}$, $m_{i}$ and $\tau_{2}$ in
(\ref{necessary_safety_active_3arcs_free}). Thus, the construction of an
optimal trajectory is obtained by solving two sub-problems and piecing the
solutions together. This is an important property because it also allows us to
easily check for the existence of a feasible solution: if $\tau_{2}<\tau_{1}$
then no feasible optimal trajectory exists in this case.
\end{customrem}





\section{Optimal Control of CAVs in the MZ}

\label{oc:mz}

As mentioned in Section \ref{modeling_turns}, the inclusion of left and right
turns requires consideration of passenger comfort which is commonly
quantified through the \emph{jerk} $J_{i}(t)=\dot{u}_{i}(t)$
\cite{hogan1984adaptive} defined as the time derivative of acceleration. Thus,
one approach to minimize passenger discomfort in the MZ is to keep the
magnitude of the resultant force, which consists of the centripetal force and
the braking force, unchanged. Note that both the magnitude of the centripetal
force and the angle between the centripetal and braking forces do not change
while the vehicle makes a turn. Hence, the following optimization problem is
formulated with the objective of minimizing the $L^{2}$-norm of jerk for each
CAV $i$:
\begin{equation}
\begin{aligned} \min_{J_{i}} & \frac{1}{2}\int_{t_{i}^{m}}^{t_{i}^{f}} J_{i}^{2}(t)dt\\ 
\text{s.t.}: &~\eqref{eq:model2}, J_{i}(t) = \dot{u}_{i}(t),\\ 
\text{given } & t_{i}^{m}, t_{i}^{f},  u_{i}(t_{i}^{m}),u_{i}(t_{i}^{f}), \\
  & v_{i}^{m}, v_{i}^{f}, p_{i}(t_{i}^{m}), p_{i}(t_{i}^{f}). \label{eq:comfort} \end{aligned}
\end{equation}
and, following the definitions in Section \ref{modeling_turns}, we have
\begin{equation}
\begin{aligned} &p_{i}(t_{i}^{m}) = L,\\ &p_{i}(t_{i}^{f})=\left\{ \begin{array} [c]{ll}L + S_L, & \mbox{if $d_i = 0$},\\ L + S, & \mbox{if $d_i = 1$},\\ L + S_R, & \mbox{if $d_i = 2$}, \end{array} \right. \label{mz_ocp:boundary_position} \end{aligned}
\end{equation}
The analytical solution of problem \eqref{eq:comfort} was obtained in
\cite{ntousakis2016optimal} using Hamiltonian analysis and considering the
jerk as the control input.

\subsection{Joint Minimization of Passenger Discomfort and Energy Consumption}

\label{oc:joint} In dealing with turning CAVs in the MZ, we formulate a joint
objective expressed as a convex combination of acceleration/deceleration and
jerk as follows:
\begin{equation}
\begin{aligned}  \min_{J_{i}} & \frac{1}{2}\int_{t_i^m}^{t_i^f} [\rho_1 u_{i}^{2}(t) + \rho_2 J_{i}^{2}(t)]dt \\ 
\text{s.t.}: &~\eqref{eq:model2}, J_{i}(t) = \dot{u}_{i}(t),\\ 
\text{given } & t_{i}^{m}, t_{i}^{f}, u_{i}(t_{i}^{m}), v_{i}^{m}, \eqref{mz_ocp:boundary_position}, v_{i}^{f}.\label{eq:comfort1} \end{aligned}
\end{equation}
where $\rho_{1}=w\cdot q_{1}$, $\rho_{2}=(1-w)\cdot q_{2}$ with $q_{1}$,
$q_{2}$ being the normalization factors which are selected so that $q_{1}\cdot
u_{i}^{2}\in\lbrack0,1]$ and $q_{2}\cdot J_{i}^{2}\in\lbrack0,1]$, and
$w\in\lbrack0,1]$ is a weight associated with the importance of energy
consumption relative to passenger discomfort. Note that $t_{i}^{m}$ has
already been determined by solving the CZ optimal control problem in the
previous section, hence also $u_{i}(t_{i}^{m})$ and $v_{i}^{m}$ are known.
Finally, $v_{i}^{f}$ is also  set to a constant $v_{i}^{f}=v^{f}$ for all
CAVs, where $v^{f}$ is a predetermined desired exit speed (e.g., $v_{\max}$ if
we wish to maximize traffic throughput after the intersection). The reason for
selecting a common speed is to prevent the chance of collisions that might
result when CAVs exit the MZ at different speeds.

Given the objective function in \eqref{eq:comfort1}, the Hamiltonian is
\begin{equation}
\begin{aligned} H_{i}(p_{i}, v_i, u_i, J_{i}, \lambda_{i}, t) &= \frac{1}{2}[\rho_{1} u_{i}^{2}(t) + \rho_{2} J_{i}^{2}(t)] \\ & + \lambda_{i}^p v_i(t) + \lambda_i^v u_i(t) + \lambda_i^u J_i(t), \label{mz:Hamiltonian}\end{aligned}
\end{equation}
where $\lambda_{i}^{p}$, $\lambda_{i}^{v}$, $\lambda_{i}^{u}$ are the costate
variables. Note that the collision avoidance inside the MZ has been ensured
through \eqref{Theorem1}. The Euler-Lagrange equations associated with the
position and speed respectively are
\begin{equation}
\dot{\lambda}_{i}^{p}=-\frac{\partial H_{i}}{\partial p_{i}}=0,~\text{\ }%
\dot{\lambda}_{i}^{v}=-\frac{\partial H_{i}}{\partial v_{i}}=-\lambda_{i}%
^{p},\label{mz_eq:lambda_pv}%
\end{equation}
and
\begin{equation}
\dot{\lambda}_{i}^{u}=-\frac{\partial H_{i}}{\partial u_{i}}=-\rho_{1}%
u_{i}-\lambda_{i}^{v},\label{mz_eq:lambda_u}%
\end{equation}
is associated with the acceleration. Since the terminal acceleration/deceleration
$u_{i}(t_{i}^{f})$ is not pre-specified, we also have the transversality
condition
\begin{equation}
\lambda_{i}^{u}(t_{i}^{f})=0.\label{mz_eq:costate_trans}%
\end{equation}
The necessary condition for optimality is
\begin{equation}
\frac{\partial H_{i}}{\partial J_{i}}=\rho_{2}J_{i}(t)+\lambda_{i}%
^{u}=0,\label{mz_eq:optimum}%
\end{equation}
Note that \eqref{mz_eq:lambda_pv} leads to
\begin{equation}
\lambda_{i}^{p}=a_{i},\text{ and }\lambda_{i}^{v}=-a_{i}t-b_{i},
\end{equation}
where $a_{i}$ and $b_{i}$ are constants of integration. Given
\eqref{mz_eq:lambda_u}, we have
\begin{equation}
\lambda_{i}^{u}=\frac{1}{2}a_{i}t^{2}+b_{i}t-\rho_{1}v_{i}+c_{i}%
\label{mz_eq:lambda_u2}%
\end{equation}
where $c_{i}$ is a constant. We can now obtain a complete analytical solution
of \eqref{eq:comfort1} as follows.

\begin{customthm}
The optimal trajectory for \eqref{eq:comfort1} is given by
\begin{align}
J_{i}^{*}(t)  & = \frac{a_{i}}{\rho_{1}} + e_{i} A_{1}^{3} e^{A_{1} t} +
f_{i}A_{2}^{3} e^{A_{2} t}\label{mz:ji*}\\
u_{i}^{*}(t)  & = \frac{1}{\rho_{1}}(a_{i} t + b_{i}) + e_{i} A_{1}^{2}
e^{A_{1} t}+ f_{i} A_{2}^{2} e^{A_{2} t}\label{mz:ui*}\\
v_{i}^{*}(t) &  = \frac{1}{\rho_{1}}(\frac{1}{2}a_{i} t^{2} + b_{i} t + c_{i}
+ \frac{a_{i} \rho_{2}}{\rho_{1}}) \nonumber\\
& + d_{i} A_{1} e^{A_{1} t} + f_{i} A_{2}
e^{A_{2} t}\label{mz:vi*}\\
p_{i}^{*}(t)  & = \frac{1}{\rho_{1}}(\frac{1}{6}a_{i} t^{3} + \frac{1}{2}b_{i}
t^{2} + c_{i} t + \frac{a_{i}\rho_{2}}{\rho_{1}}t + d_{i})\nonumber\\
&  + e_{i} e^{A_{1} t} + f_{i} e^{A_{2} t}\label{mz:pi*}
\end{align}
where $A_{1} = \sqrt{\frac{\rho_{1}}{\rho_{2}}}$, $A_{2} = -\sqrt{\frac
{\rho_{1}}{\rho_{2}}}$, and $a_{i}$, $b_{i}$, $c_{i}$, $d_{i}$, $e_{i}$ and
$f_{i}$ are constants of integration determined through
\begin{subequations}
\begin{align}
\frac{1}{\rho_{1}}(\frac{1}{6}a_{i} \cdot(t_{i}^{m})^{3} + \frac{1}{2}b_{i}
\cdot(t_{i}^{m})^{2} + c_{i} t_{i}^{m}  & \nonumber\\
+ \frac{a_{i}\rho_{2}}{\rho_{1}}t_{i}^{m} + d_{i}) + e_{i} e^{A_{1} t_{i}^{m}}
+ f_{i} e^{A_{2} t_{i}^{m}}  &  = L\label{mz:equations:pm}\\
\frac{1}{\rho_{1}}(\frac{1}{2}a_{i} \cdot(t_{i}^{m})^{2} + b_{i} t_{i}^{m} +
c_{i} + \frac{a_{i} \rho_{2}}{\rho_{1}})  & \nonumber\\
+ d_{i} A_{1} e^{A_{1} t_{i}^{m}} + f_{i} A_{2} e^{A_{2} t_{i}^{m}}  & =
v_{i}^{m}\label{mz:equations:vm}\\
\frac{1}{\rho_{1}}(a_{i} t_{i}^{m} + b_{i}) + e_{i} A_{1}^{2} e^{A_{1}
t_{i}^{m}}+ f_{i} A_{2}^{2} e^{A_{2} t_{i}^{m}}  &  = u_{i}^{m}%
\label{mz:equations:um}\\
\frac{1}{\rho_{1}}(\frac{1}{6}a_{i} \cdot(t_{i}^{f})^{3} + \frac{1}{2}b_{i}
\cdot(t_{i}^{f})^{2} + c_{i} t_{i}^{f}  & \nonumber\\
+ \frac{a_{i}\rho_{2}}{\rho_{1}}t_{i}^{f} + d_{i}) + e_{i} e^{A_{1} t_{i}^{f}}
+ f_{i} e^{A_{2} t_{i}^{f}}  &  = p_{i}^{f}\label{mz:equations:pf}\\
\frac{1}{\rho_{1}}(\frac{1}{2} a_{i} \cdot(t_{i}^{f})^{2} + b_{i} t_{i}^{f} +
c_{i} + \frac{a_{i} \rho_{2}}{\rho_{1}})  & \nonumber\\
+ d_{i} A_{1} e^{A_{1} t_{i}^{f}} + f_{i} A_{2} e^{A_{2} t_{i}^{f}}  &  =
v_{i}^{f}\label{mz:equations:vf}\\
\frac{a_{i}}{\rho_{1}} + e_{i} A_{1}^{3} e^{A_{1} t_{i}^{f}} + f_{i}A_{2}^{3}
e^{A_{2} t_{i}^{f}}  & = 0.\label{mz:equations:vf_free}%
\end{align}
\label{mz:equations}
\end{subequations}
\end{customthm}

\emph{Proof. } Combining \eqref{mz_eq:lambda_u2} with \eqref{mz_eq:optimum},
we have the ordinary differential equation:
\begin{equation}
\rho_{2}\ddot{v_{i}}-\rho_{1}v_{i}+\frac{1}{2}a_{i}t^{2}+b_{i}t+c_{i}%
=0.\nonumber
\end{equation}
whose solution yields the optimal speed \eqref{mz:vi*}. Using \eqref{mz:vi*}
in the system dynamics, we can then derive \eqref{mz:ji*}, \eqref{mz:ui*} and
\eqref{mz:pi*}. The first five equations in \eqref{mz:equations} follow from
the boundary conditions $p_{i}(t_{i}^{m})=L$, $v_{i}(t_{i}^{m})$ and
$u_{i}(t_{i}^{m})$ (known by $t_{i}^{m}$), and the specified $p_{i}(t_{i}%
^{f})$ and $v_{i}(t_{i}^{f})$. The last equation follows from the
transversality condition \eqref{mz_eq:costate_trans}:
\begin{align*}
&  \lambda_{i}^{u}(t_{i}^{f})=-\frac{J_{i}(t_{i}^{f})}{\rho_{2}}\\
&  =\rho_{2}[\frac{a_{i}}{\rho_{1}}+e_{i}A_{1}^{3}e^{A_{1}t_{i}^{f}}%
+f_{i}A_{2}^{3}e^{A_{2}t_{i}^{f}}]=0.
\end{align*}
using \eqref{mz:ji*}. \hfill$\blacksquare$

Note that since $0\leq w\leq1$, the optimal solution is only valid when
$w\neq1$ and $w\neq0$. When $w=0$, the problem becomes \eqref{eq:comfort} with
the objective of minimizing jerk only. When $w=1$, the problem minimizes
energy consumption only. Although the state and control constraints are not
incorporated in \eqref{eq:comfort1}, it is possible that the minimum speed $v_{min}$ and/or maximum deceleration $u_{min}$ constraints become active. The approach to analyze such cases is similar to the analysis in Appendix \ref{appendix:CZ}.


\section{Simulation Examples\label{analytic}}

We begin with several numerical examples illustrating the different cases
discussed in Section \ref{CZ}. In terms of computational complexity, we should
point out that except for the case where the complete solution is given by the
simple system of linear equations (\ref{fixed_term_time_solution}), solving a
system of nonlinear equations with six parameters involved as in
(\ref{necessary_safety_active_3arcs_free}) is certainly nontrivial. A good
initial `guess' of the parameter values is extremely useful in the convergence
of the root-finding algorithm for numerical solvers. To do so, our approach is
to start with the unconstrained optimal solution and add the constraints step
by step. At each step, when a constraint is added to the problem, we use the
solution from the previous step as the initial estimate and obtain the new
solution, which is then used as the initial estimate for the problem
formulated in the next step; if the current solution satisfies all the
constraints, then we have obtained the optimal solution.

\begin{figure}[ptb]
\centering
\includegraphics[width= 1\columnwidth]{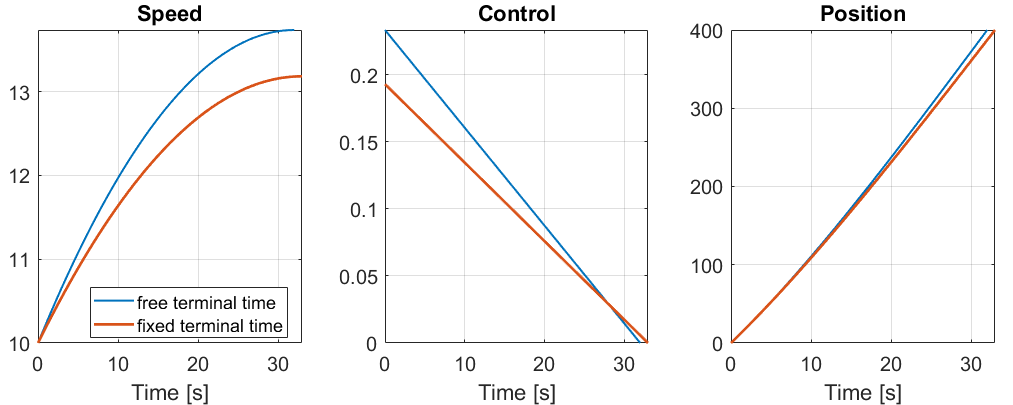}
\caption{Unconstrained optimal trajectories with free and fixed terminal
times.}%
\label{no_constraint}%
\end{figure}

\textit{Unconstrained optimal control with free terminal time.} The parameters
used are: $L=400$m, $\gamma=0.1$, $v_{i}^{0}=10$m/s, $t_{i}^{0}=0$s. The
optimal terminal time is obtained as $t_{i}^{m}=32.03$s as shown by the blue
curves in Fig. \ref{no_constraint}.

\textit{Unconstrained optimal control with fixed terminal time.} Assuming
$t_{i}^{m}=32.03$s violates \eqref{Theorem1}, and we need to formulate
problem $P_{1}$ by adding $t_{i}^{m}=33$s to $P_{0}$. The resulting optimal
control, speed, and position trajectories are shown by the red curves in Fig.
\ref{no_constraint}.

\begin{figure}[ptb]
\centering
\includegraphics[width= 1\columnwidth]{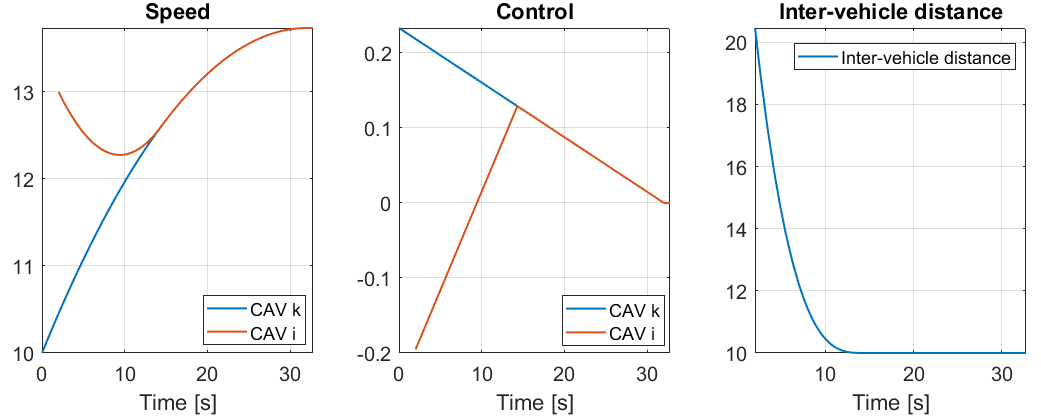}\caption{The
state constraint $p_{i}(t)+\delta-p_{k}(t)\leq0$ active (no exit).}%
\label{safety:simulation_2arcs}%
\end{figure}

\textit{Safety-constrained optimal control without exit.}
Assuming CAV $k=1$ enters the CZ at $t_{k}^{0}=0$ with an initial speed
$v_{k}^{0}=10$m/s and exits at $t_k^m=32.03$s, the optimal control is $u_{k}%
^{\ast}(t)=-0.0073t+0.23$. Then, we assume that CAV $i=2$ enters the CZ at
$t_{i}^{0}=2$s with an initial speed $v_{i}^{0}=13$m/s. The terminal time of
CAV $i$ is $t_{i}^{m}=t_{k}^{m}+\frac{\delta}{v_{k}^{m}}=32.76$s where the
minimal safe following distance is $\delta=10$m. The optimal control for CAV
$i$ is
\begin{equation}
u_{i}^{\ast}(t)=\left\{
\begin{array}
[c]{ll}%
0.0263t-0.25, & \mbox{$t \in [2, 14.31]$}\\
-0.0073t+0.23, & \mbox{$t \in (14.31, 32.03]$},\\
0, & \mbox{$t \in (32.03, 32.76]$},
\end{array}
\right.  \nonumber
\end{equation}
as shown in Fig. \ref{safety:simulation_2arcs}. Note that in this case, CAV
$i$ needs to start out by decelerating before entering the constrained arc.

\begin{figure}[ptb]
\centering
\includegraphics[width= 1\columnwidth]{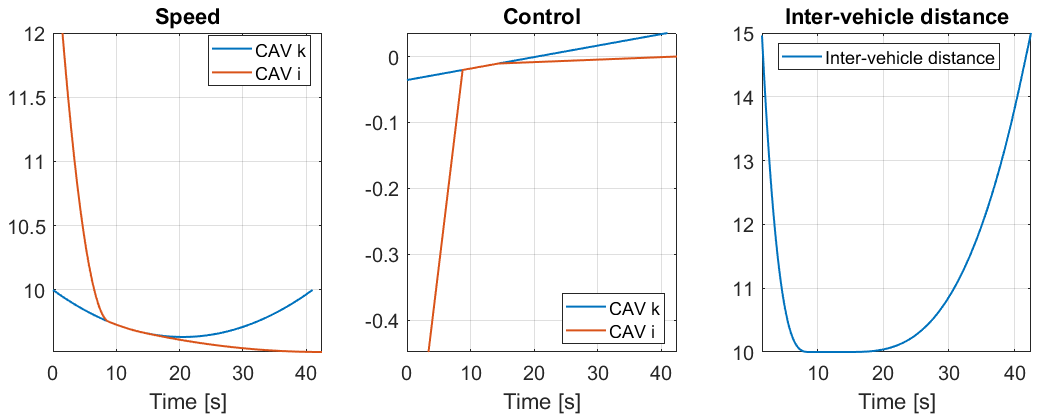}\caption{The
state constraint $p_{i}(t) + \delta- p_{k}(t) \leq0$ active (with entry and
exit).}%
\label{safety:simulation_3arcs}%
\end{figure}

\textit{Safety-constrained optimal control with exit.} Assuming
CAV $k=1$ enters the CZ at $t_{k}^{0}=0$ with an initial speed $v_{k}^{0}=10$
and exits at $t_{k}^{m}=41$s with a terminal speed $v_{k}^{m}=10$m/s, the
optimal control is $u_{k}^{\ast}(t)=0.0017t-0.0357$. Then, we assume that CAV
$i=2$ enters the CZ at $t_{i}^{0}=1.5$s with an initial speed $v_{i}^{0}%
=12$m/s, and the terminal time of CAV $i$ is $t_{i}^{m}=42.5$s satisfying
$t_{i}^{m}>t_{k}^{m}+\delta/v_{k}^{m}$ where the minimal safe following
distance is $\delta=10$m. The optimal control for CAV $i$ is
\begin{equation}
u_{i}^{\ast}(t)=\left\{
\begin{array}
[c]{ll}%
0.07971t-0.7183, & \mbox{$t \in [1.5, 8.75]$}\\
0.0017t-0.0357, & \mbox{$t \in (8.75, 14.4]$},\\
0.00038t-0.0161 & \mbox{$t \in (14.4, 42.5]$},
\end{array}
\right.  \nonumber
\end{equation}
as shown in Fig. \ref{safety:simulation_3arcs}. 

\subsection{Optimal Trajectories in the CZ \label{sim_CZ}}

\begin{figure}[ptb]
\centering
\includegraphics[width= 1\columnwidth]{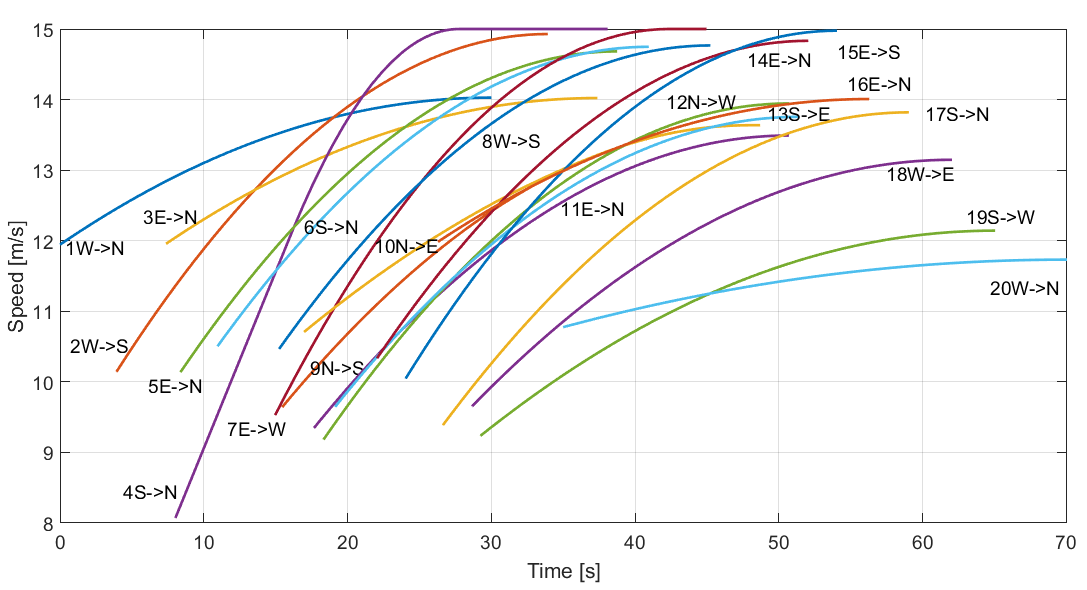} \caption{Speed profiles
of the first 20 CAVs in the CZ.}%
\label{vit_cz}%
\end{figure}

\begin{figure}[ptb]
\centering
\includegraphics[width= 1\columnwidth]{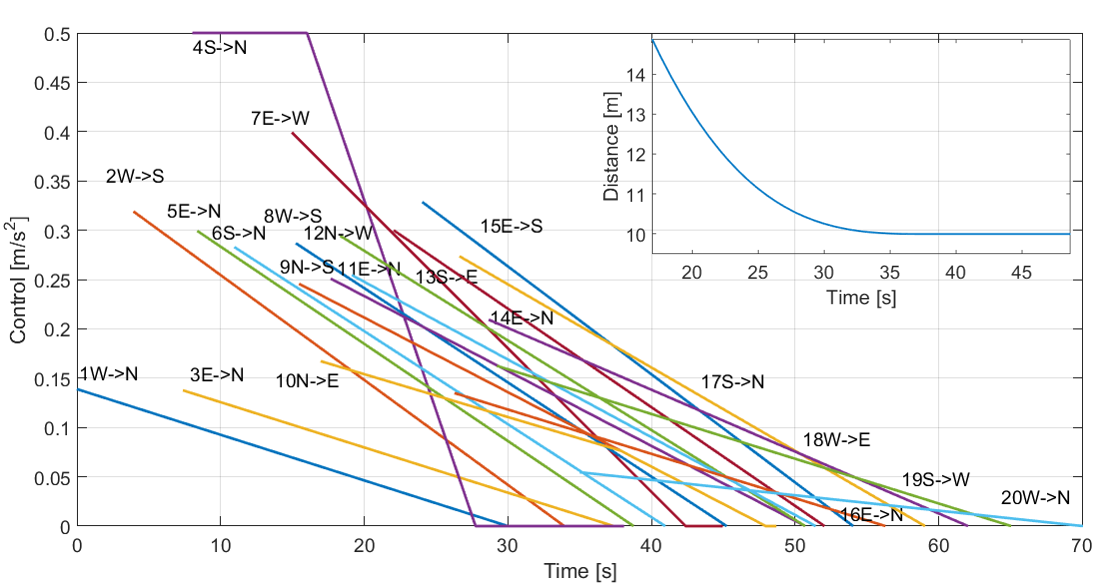} \caption{Control
input/acceleration profiles of the first 20 CAVs in the CZ.}%
\label{uit_cz}%
\end{figure}

The proposed decentralized optimal control framework inside the CZ is
illustrated through simulation in MATLAB. We assumed a single lane for each
traffic direction and the parameters used are: $L=400$m, $S=30$m; the speed
constraints are $v_{max}=15$m/s and $v_{min}=5$m/s; the control constraints
are $u_{max}=0.5$m/s$^{2}$ and $u_{min}=-0.5$m/s$^{2}$; $S_{L}=\frac{3}{8}\pi
S$, $S_{R}=\frac{1}{8}\pi S$, $\delta=10$ m, and $\Delta_{i}=5,3,3$s for a
left turn, going straight, and a right turn respectively. CAVs arrive at the
CZ based on a random arrival process which we assumed to be a Poisson process
with rate $\lambda=1$ and the initial speeds are uniformly distributed over
$[8,12]$m/s.

The optimal speed and control trajectories in the CZ are shown in Figs.
\ref{vit_cz} and \ref{uit_cz}, with labels indicating the position of the CAV
in the FIFO queue and the driving direction. To illustrate cases where the
acceleration and speed constraints become active, note that for CAV \#4, the
optimal trajectory consists of three arcs: $(i)$ a constrained arc starting at
$t_{4}^{0}=8.03$s where the CAV accelerates at $u_{max}$ until $\tau
_{1}=16.02$s, $(ii)$ an unconstrained arc where the CAV accelerates to
$v_{max}$ at $\tau_{2}=27.57$s, $(iii)$ a constrained arc where the CAV
cruises at $v_{max}$ ($u_{i}(t)=0$) until it exits the CZ at $t_{4}^{m}%
=38.09$s. For CAV \#7, the optimal trajectory consists of two arcs: $(i)$ an
unconstrained arc starting at $t_{7}^{0}=14.96$s where the CAV keeps
accelerating until it reaches $v_{max}$ at $\tau=42.38$s, $(ii)$ a constrained
arc where the CAV cruises at $v_{max}$ ($u_{i}(t)=0$) until $t_{7}^{m}=44.96$s.

To illustrate a case where the safety constraint is included, observe that CAV
\#10 is traveling on the same lane as \#9. At $\tau=37.85$s, \#10 enters the
constrained arc where $p_{10}(t)+10-p_{9}(t)=0$. It then follows the optimal
trajectory of \#9, which can be reflected by the slope change in Fig.
\ref{uit_cz}. The inter-vehicle distance between \#9 and \#10 is shown as a
subfigure in Fig. \ref{uit_cz}. Observe that after \#10 enters the constrained
arc, it stays constrained until reaching the MZ without exiting. The lower
bound of $t_{10}^{m}$ happens to be constrained by \#9, hence, once \#10
enters the constrained arc, there is no incentive for it to exit.



\begin{figure}[ptb]
\centering
\includegraphics[width= 0.8\columnwidth]{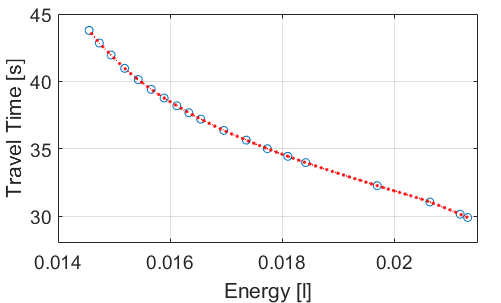} \caption{Pareto
efficiency sets and frontier corresponding to different combinations of energy
consumption and traffic throughput in the CZ.}%
\label{pareto_cz}%
\end{figure}

The optimal solution to \eqref{eq:p_ori} varies as the weight $\beta$ changes.
To investigate the tradeoff between energy consumption and traffic throughput, we examine a range of cases with different $\beta$ values and generate the associated Pareto sets illustrating the fact that no objective can be made
better off without making at least one other objective worse. In \eqref{eq:p_ori}, we use $u_i^2(t)$ as a rough approximation of energy consumption, since it adequately captures its monotonic dependence on acceleration, while also allowing us to derive the analytical solution. However, to more accurately assess the impact of our controller, we use the polynomial metamodel proposed in \cite{Kamal2013a} which yields vehicle energy consumption as a function of speed and acceleration: $\dot{f}=\dot{f}_{cruise}+\dot{f}_{accel}$ where $\dot{f}_{cruise}=w_{0}+w_{1}v_{i}(t)+w_{2}v_{i}^{2}(t)+w_{3}v_{i}^{3}(t)$ estimates the energy consumed by a vehicle cruising at a constant
speed $v_{i}(t)$, and $\dot{f}_{accel}=u_{i}(t)\cdot\lbrack r_{0}+r_{1}%
v_{i}(t)+r_{2}v_{i}^{2}(t)]$ estimates the additional energy consumed due to
acceleration with $u_{i}(t)$. The polynomial coefficients $w=[w_{0}%
,w_{1},w_{2},w_{3}]$ and $r=[r_{0},r_{1},r_{2}]$ are calculated from
experimental data. In addition, we use the average travel time inside the CZ, i.e., $t_i^m - t_i^0$ as a measurement of the traffic throughput.
By obtaining all
of the optimal solutions to \eqref{eq:p_ori} while varying the weight $\beta$,
we can derive the Pareto sets and the Pareto frontier corresponding to
different combinations of fuel consumption and average travel time as shown
in Fig. \ref{pareto_cz}. Observe that there exists a tradeoff between minimizing energy consumption and maximizing traffic throughput. 

To evaluate the effectiveness of the proposed solution, we carried out a comparison with the baseline scenario using VISSIM, where all the vehicles are assumed to be non-CAVs under the control of traffic lights with fixed switching times.
The comparison is shown in Table \ref{comp_baseline}, where the weight $\beta$ in \eqref{eq:decentral}, used for trading off energy and throughput, is set to 0.75, 0.5, and 0.25 respectively. When $\beta=0.5$ where energy and throughput are equally weighted, the energy consumption improvement is 13.46\%, while the average travel time is improved by 29.84\% compared to the baseline scenario. 
As $\beta$ varies (see Table \ref{comp_baseline}) the resulting tradeoff between travel time and energy changes as expected. Since our objective is to jointly minimize energy consumption and maximize traffic throughput, we also compute the value of the objective function in \eqref{eq:decentral}. As shown in Table \ref{comp_baseline}, the optimized non-signalized performance is significantly better than the signalized baseline regardless of $\beta$ values.


\begin{table}[ptb]
\caption{Comparison with baseline scenario (signalized intersection)}%
\label{comp_baseline}
\begin{tabularx}{\linewidth}{l l l l l}
\toprule  & OC\_${\beta_1}$ & OC\_${\beta_2}$  & OC\_${\beta_3}$ & Baseline \\
\midrule Travel Time [s] &  28.92   & 30.99 &  34.98   & 44.17 \\
Energy [l] &  0.047 & 0.045 &  0.042   & 0.052 \\
OC Cost\_${\beta_1}$ &  11.42  & / & /  & 16.57 \\
OC Cost\_${\beta_2}$ & / & 4.08 & /  & 5.52 \\
OC Cost\_${\beta_3}$  & / & / &  1.50  & 1.84 \\
\bottomrule
\end{tabularx}
\footnote{1} $\beta_1 = 0.75$,  $\beta_2=0.5$, $\beta_3=0.25$.
\end{table}

\subsection{Optimal Trajectories in the MZ \label{sim_MZ}}



\begin{figure}[tbh]
\centering
\includegraphics[width= 1\columnwidth]{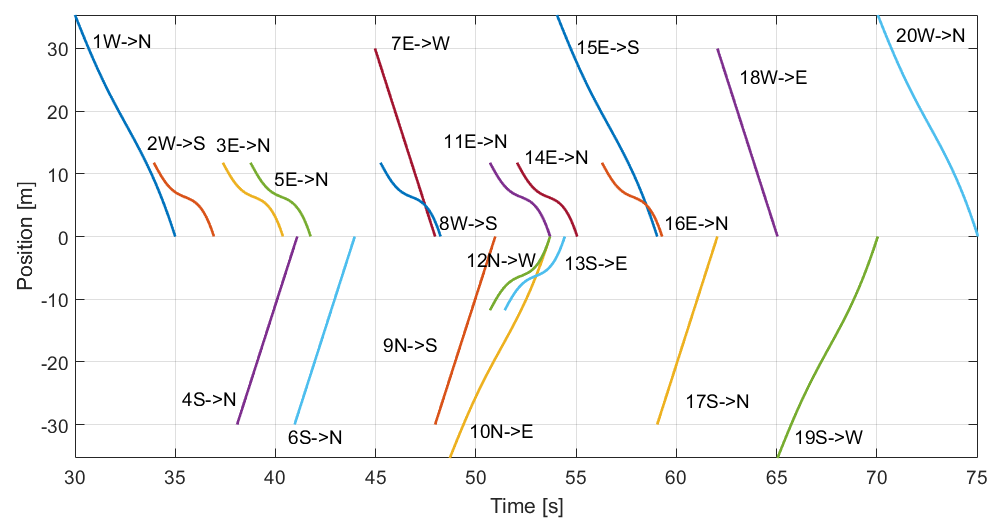} \caption{Distance to the
end of MZ of the first 20 CAVs in the MZ.}%
\label{p_MZ}%
\end{figure}


The proposed decentralized optimal control framework inside the MZ
incorporating turns is illustrated through simulation in MATLAB with the same
model parameters as in Section \ref{sim_CZ}. For simulation purposes, we
assume the speeds at the entry of the MZ are set to $8$m/s for CAVs turning
left, $6$m/s for CAVs turning right, and $10$m/s for CAVs going straight,
respectively. The speed at the exit of the MZ is set to $v^f = 10$m/s.

The position trajectories of the first 20 CAVs inside the MZ are shown in Fig.
\ref{p_MZ}. CAVs are separated into two groups: those shown above zero are
driving from east or west, and those below zero are driving from north or
south, with labels indicating the position of the vehicles in the FIFO queue
and the driving direction. Observe that CAV \#11 belongs to $\mathcal{O}%
_{12}(t)$ and no collision would occur between \#11 and \#12. Hence, they can be
traveling inside the MZ at the same time, i.e., $t_{12}^{f}=t_{11}^{f}$.
Similarly, since CAV \#12 belongs to $\mathcal{O}_{13}(t)$, no collision would
occur between \#12 and \#13 as well. However, recalling the dependence of the
terminal conditions on $e$, $s$, $l$, $o$ in \eqref{Theorem1}, \#13 is
constrained by \#10 which may collide with it at the end of the MZ. As a
result, \#13 has to wait until \#10 leaves the MZ for a distance $\delta$,
which leads to $t_{13}^{f}=t_{10}^{f}+\frac{\delta}{v_{10}^{f}}$.

\begin{figure}[ptb]
\centering
\includegraphics[width= 0.8\columnwidth]{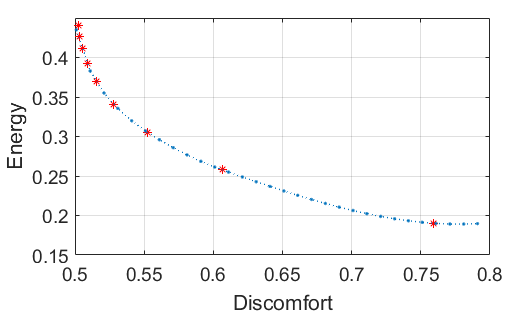} \caption{Pareto
efficiency sets and frontier corresponding to different combinations of energy
consumption and passenger discomfort in the MZ.}%
\label{pareto}%
\end{figure}

The optimal solution to \eqref{eq:comfort1} varies as the weight $w$ changes.
Similarly, to illustrate the tradeoff between passenger discomfort and energy
consumption, we can examine a range of cases with different weights and generate
the associated Pareto sets. In Fig. \ref{pareto}, we use $u_i^2(t)$ and $J_i^2(t)$ as rough approximations of energy consumption and passenger discomfort respectively. Clearly, there is a tradeoff between these two metrics
as the vehicle turns within the intersection.

\section{Conclusions and Future Work}

We have extended earlier work in \cite{ZhangMalikopoulosCassandras2016} and
\cite{malikopoulos2018decentralized} where a decentralized optimal control framework was
established  for optimally controlling CAVs crossing an urban intersection
without considering turns. In this paper, we have included left and right
turns and considered a new optimal control problem formulation where the
tradeoff between energy and travel time is explicitly quantified for CAVs in
the CZ and all safety constraints are incorporated. Despite the added
complexity of turns, we have shown that the optimal solution can still be
obtained in decentralized fashion, with each CAV requiring information from a
subset of other CAVs. This enables the on-line solution to be obtained by
on-board computation resources for each individual CAV. In addition, we have
formulated another optimization problem with the objective of minimizing
passenger discomfort while the vehicle turns at the MZ of the intersection,
and investigated the tradeoff between minimizing energy consumption and
passenger discomfort.
Ongoing research is exploring the effect of partial CAV penetration in mixed
traffic situations where both CAVs and human-driven vehicles share the the
road \cite{Zhang2018penetration}. Future work will investigate the effect of
errors and/or communication delays to the system behavior as well as the
coupling between multiple intersections, and the possibility of extending the
resequencing approach in \cite{Zhang2018sequence} to potentially improve
overall traffic throughput.


\useRomanappendicesfalse
\appendices

\section{Constrained Optimal Control Analysis \label{appendix:CZ}}

In this Appendix, we discuss the effect of the control and state constraints
not included in our analysis of Section \ref{CZ_solution}. Clearly, there are
a number of possible situations that may arise, including the possibility of
both the state constraint $v_{i}(t)-v_{max}\leq0$ and the control constraint
$u_{i}(t)-u_{max}\leq0$ becoming active (e.g., CAV \#4 in Figs. \ref{vit_cz}
and \ref{uit_cz}). The analysis provided here is limited to the state
constraint $v_{i}(t)-v_{max}\leq0$ or the control constraint $u_{i}%
(t)-u_{max}\leq0$ each becoming active on its own. These basic cases serve as
\textquotedblleft building blocks\textquotedblright\ to handle situations of
multiple constraints becoming active when this is feasible.

\subsection{State constraint $v_{i}(t) - v_{max} \leq 0$ active only}

\label{vmax_active}

When the state constraint $v_{i}(t)-v_{max}\leq0$ becomes active at $\tau$,
the jump conditions in \eqref{ipm_disc} become
\begin{equation}
\begin{aligned} \lambda_i^{p}(\tau^{-}) &= \lambda_i^{p}(\tau^{+}) +
\pi_i^p \frac{\partial (v_i(\tau)-v_{max})}{\partial p_i} = \lambda_i^{p}(\tau^{+}),\\
\lambda_i^{v}(\tau^{-}) &= \lambda_i^{v}(\tau^{+}) + \pi_i^v \frac{\partial (v_i(\tau)-v_{max})}{\partial v_i}\\
& = \lambda_i^{v}(\tau^{+}) + \pi_i^v(\tau),\\
H_i(\tau^{-}) &= H_i(\tau^{+}) - \pi_i(v_i(\tau)-v_{max})_{t} = H_i(\tau^{+}),\\ & \pi_i(v_i(\tau)-v_{max}) = 0, \pi_i \geq0. \end{aligned}\label{corner}%
\end{equation}
Thus, we have $H_{i}(\tau^{-})=H(\tau^{+})$ and $\lambda_{i}^{p}(\tau
^{-})=\lambda_{i}^{p}(\tau^{+})$. Based on this fact, we have the following
result whose proof is similar to that of Theorem \ref{safety:u_cont}.

\begin{customthm}
The optimal control $u_{i}(t)$ is continuous at the junction $\tau$ of the unconstrained and $v_{max}$-constrained arcs, i.e., 
$u_i^*(\tau^-) = u_i^*(\tau^+)$. \label{vmax_active:u_continuous}
\end{customthm}

For this case, CAV $i$ remains on the constrained arc until it reaches the MZ
and we have
\begin{equation}
u_{i}(t)=\left\{
\begin{array}
[c]{ll}%
a_{i}t+b_{i}, & \mbox{$t \leq \tau $},\\
0, & \mbox{$t > \tau$}.
\end{array}
\right.  \label{u_2arcs}%
\end{equation}
where $\tau$ is the entry point of the constrained arc $v_{i}(t)-v_{max}=0$,
and, due to Theorem \ref{vmax_active:u_continuous}, we also have $a_{i}%
\tau+b_{i}=0$.

Combined with the boundary conditions and the transversality conditions
\eqref{eq:costate_trans}, we have the following conditions
\begin{subequations}
\begin{align}
\frac{1}{6}a_{i}\cdot(t_{i}^{0})^{3}+\frac{1}{2}b_{i}\cdot(t_{i}^{0}%
)^{2}+c_{i}t_{i}^{0}+d_{i} &  =0\label{necessary_vmax_active:p0}\\
\frac{1}{2}a_{i}\cdot(t_{i}^{0})^{2}+b_{i}t_{i}^{0}+c_{i} &  =v_{i}%
^{0}\label{necessary_vmax_active:v0}\\
\frac{1}{6}a_{i}\tau^{3}+\frac{1}{2}b_{i}\tau^{2}+c_{i}\tau+d_{i}-L &
~~\nonumber\\
-v_{max}(t_{i}^{m}-\tau) &  =0\label{necessary_vmax_active:p_tau}\\
\frac{1}{2}a_{i}\tau^{2}+b_{i}\tau+c_{i} &  =v_{max}%
\label{necessary_vmax_active:v_tau}\\
a_{i}\tau+b_{i} &  =0\label{necessary_vmax_active:u_tau}\\
\gamma+a_{i}\cdot v_{max} &  =0\label{necessary_vmax_active:Ham}%
\end{align}
\label{necessary_vmax_active}%
\end{subequations}
where $a_{i}$, $b_{i}$, $c_{i}$, $d_{i}$, $\tau$ and $t_{i}^{m}$ are obtained
by solving these six equations. The first two immediately follow from the
initial conditions, while \eqref{necessary_vmax_active:p_tau} follows from
$p_{i}(\tau)+v_{max}\cdot(t_{i}^{m}-\tau)=L$,
\eqref{necessary_vmax_active:v_tau} follows from $v_{i}(\tau)=v_{max}$, and
\eqref{necessary_vmax_active:u_tau} follows from the fact that $u_{i}(\tau
^{-})=u_{i}(\tau^{+})=0$. The last equation \eqref{necessary_vmax_active:Ham}
follows from the continuity of the Hamiltonian and the transversality
condition $H_{i}(t_{i}^{m})=0$. Note that in this case we assumed a free
terminal time $t_{i}^{m}$. If $t_{i}^{m}$ is fixed, then we simply remove the
transversality condition \eqref{necessary_vmax_active:Ham}. Simulation
example of this case arises in Figs. \ref{vit_cz} and \ref{uit_cz}, e.g., CAV \#7.

\subsection{Control constraint $u_{i}(t) - u_{max} \leq0$ active only}

\label{umax_active} When the control constraint $u_{i}(t)-u_{max}\leq0$ is
active at $\tau$, both the Hamiltonian and the costates are continuous
according to \eqref{ipm_disc}. For this case, we have
\begin{equation}
u_{i}(t)=\left\{
\begin{array}
[c]{ll}%
u_{max}, & \mbox{$ t \leq \tau $},\\
a_{i}t+b_{i}, & \mbox{$t > \tau$}.
\end{array}
\right.
\end{equation}
Combined with the boundary conditions and the transversality conditions
\eqref{eq:costate_trans}, we have the following conditions
\begin{subequations}
\begin{align}
(\tau-t_{i}^{0})[v_{i}^{0}+\frac{1}{2}u_{max}(\tau-t_{i}^{0})]~~~~~ &
\nonumber\\
-(\frac{1}{6}a_{i}\tau^{3}+\frac{1}{2}b_{i}\tau^{2}+c_{i}\tau+d_{i}) &
=0\label{necessary_umax_active:p_tau}\\
v_{i}^{0}+u_{max}(\tau-t_{i}^{0})-(\frac{1}{2}a_{i}\tau^{2}+b_{i}\tau+c_{i})
&  =0\label{necessary_umax_active:v_tau}\\
a_{i}\tau+b_{i} &  =u_{max}\label{necessary_umax_active:u_tau}\\
\frac{1}{6}a_{i}\cdot(t_{i}^{m})^{3}+\frac{1}{2}b_{i}\cdot(t_{i}^{m}%
)^{2}+c_{i}t_{i}^{m}+d_{i} &  =L\label{necessary_umax_active:pm}\\
a_{i}t_{i}^{m}+b_{i} &  =0\label{necessary_umax_active:vm_free}\\
\gamma- \frac{1}{2}b_{i}^{2} + a_{i} c_{i}  &  = 0\label{necessary_umax_active:Ham}%
\end{align}
\label{necessary_umax_active}%
\end{subequations}
where $a_{i}$, $b_{i}$, $c_{i}$, $d_{i}$, $\tau$ and $t_{i}^{m}$ are obtained
by solving these equations. The first two equations \eqref{necessary_umax_active:p_tau} and \eqref{necessary_umax_active:v_tau} follow from the fact that $p_i(\tau^-) = p_i(\tau^+)$ and  $v_i(\tau^-) = v_i(\tau^+)$, \eqref{necessary_umax_active:u_tau}
follows from the fact that $u_{i}(\tau^{-})=u_{i}(\tau^{+})=u_{max}$,
\eqref{necessary_umax_active:pm} follows from the terminal condition, and the
last two equations follow from the transversality conditions. In this case, we assume a free terminal time $t_{i}^{m}$. If
$t_{i}^{m}$ is fixed, then we simply remove the last transversality condition \eqref{necessary_umax_active:Ham}.

\bibliographystyle{IEEETran}
\bibliography{LRturns}

\end{document}